%% file: arXiv.tex
\documentclass[a4paper,10pt,twoside]{amsart}
\usepackage{a4wide}
\usepackage{amsmath}
\usepackage{amsthm}
\usepackage{amsfonts}
\usepackage{amssymb}
\usepackage{mathrsfs}
\usepackage{tikz}
\usetikzlibrary{arrows,positioning,decorations.text,calc}
\usepackage{url}
\usepackage{graphicx}
\usepackage[labelfont=up]{caption}
\usepackage{subcaption} 
\usepackage{todonotes} 
\usepackage[hidelinks]{hyperref} 

\theoremstyle{plain}
\newtheorem{theorem}{Theorem}[section]
\newtheorem*{lemma*}{Lemma}
\newtheorem{lemma}[theorem]{Lemma}
\newtheorem{corollary}[theorem]{Corollary}
\newtheorem{proposition}[theorem]{Proposition}

\theoremstyle{definition}
\newtheorem{definition}[theorem]{Definition}
\newtheorem{remark}[theorem]{Remark}

\def \B {\mathcal B}
\def \C {\mathscr C}
\def \G {\mathscr G}
\def \x {\mathbf x}
\def \y {\mathbf y}
\def \proj #1{\widehat{#1}}
\def \normal #1{\overline{#1}}
\def \pack {\mathscr P}

\definecolor{darkblue}{rgb}{0,0,0.7} 
\newcommand{\darkblue}{\color{darkblue}} 
\newcommand{\Dfn}[1]{\emph{\darkblue #1}} 

\DeclareMathOperator{\cone}{\mathsf{cone}}
\DeclareMathOperator{\affine}{\mathsf{aff}}
\DeclareMathOperator{\convex}{\mathsf{conv}}
\DeclareMathOperator{\depth}{\mathsf{dp}}
\DeclareMathOperator{\prefix}{\mathsf{pre}}
\DeclareMathOperator{\ball}{\mathsf{Ball}}

\title[Boyd--Maxwell packings]{Lorentzian Coxeter Systems and Boyd--Maxwell Ball Packings}

\author[H.~Chen]{Hao Chen}
\address[H.~Chen]{Freie Universit\"at Berlin, Institut f\"ur Mathematik, Arnimallee 2, 14195 Berlin, Deutschland}
\email{haochen@math.fu-berlin.de}
\urladdr{http://page.mi.fu-berlin.de/haochen}
\thanks{H.~Chen is supported by the Deutsche Forschungsgemeinschaft within the Research Training Group ``Methods for Discrete Structures'' (GRK 1408).  An alternative version of this paper appeared in the PhD thesis of the first author~\cite{chen_thesis_2014}.}

\author[J.-P.~Labb\'e]{Jean-Philippe Labb\'e}
\address[J.-P. Labb\'e]{Freie Universit\"at Berlin, Institut f\"ur Mathematik, Arnimallee~2, 14195 Berlin, Deutschland}
\email{labbe@math.fu-berlin.de}
\urladdr{http://page.mi.fu-berlin.de/labbe}
\thanks{J.-P.~Labb\'e is supported by a FQRNT Doctoral scholarship and SFB Transregio ``Discretization in Geometry and Dynamics'' (TRR 109).}


\begin{document}

\begin{abstract}
  In the recent study of infinite root systems, fractal patterns of ball packings were observed while visualizing roots in affine space.  In this paper, we show that the observed fractals are exactly the ball packings described by Boyd and Maxwell. This correspondence is a corollary of a more fundamental result: Given a geometric representation of a Coxeter group in a Lorentz space, the set of limit directions of weights equals the set of limit roots. Additionally, we use Coxeter complexes to describe tangency graphs of the corresponding Boyd--Maxwell ball packings. Finally, we enumerate all the Coxeter systems that generate Boyd--Maxwell ball packings. 
\end{abstract}

\maketitle

\section{Introduction}\label{sec:Intro}
\input{Intro}

\section{Coxeter groups, limit roots and Boyd--Maxwell Packings}\label{sec:Prelim}
\input{Prelim}

\section{Relation between limit roots and Boyd--Maxwell Packings}\label{sec:Main}
\input{Main}

\section{Enumeration of Coxeter graphs of level 2}\label{sec:List}
\input{List}

\section*{Acknowledgement}

The authors are grateful to George Maxwell for his great availability to check the enumeration results. We also thank Christian Stump for helpful discussions, and Christophe Hohlweg and Vivien Ripoll for helpful comments on a preliminary version of this manuscript.

\bibliographystyle{alpha}
\bibliography{biblio}

\section*{Erratum to Section~\ref{ssec:cluster}}
\input{Erratum}

\end{document}

%% file: Intro.tex
We establish a connection between two seemingly unrelated concepts: infinite root systems in the Lorentz space, and a special class of ball packings initially studied by Boyd and Maxwell, which generalizes Apollonian ball packings.

A Coxeter group is usually represented as a reflection group acting on a vector space, which allows us to associate a root system to the Coxeter system; see~\cite{bourbaki_elements_1968} and~\cite{humphreys_reflection_1992}.  For infinite Coxeter systems, Vinberg introduced a more flexible geometric representation that depends on a bilinear form associated to the Coxeter system~\cite{vinberg_discrete_1971, krammer_conjugacy_2009}.  In this framework, limit roots are the accumulation points of the directions of the roots. The notion was introduced and studied in~\cite{hohlweg_asymptotical_2013}.  Properties of limit roots of infinite Coxeter systems were investigated in a series of papers.  Limit roots lie on the isotropic cone of the bilinear form associated to the geometric representation~\cite{hohlweg_asymptotical_2013}.  The cone over limit roots is the imaginary cone~\cite{dyer_imaginary_2013}.  The relations between limit roots and the imaginary cone are further investigated in~\cite{dyer_imaginary2_2013}.  

\begin{figure}[!ht]
  \centering
  \begin{tabular}{p{.45\linewidth}@{\hspace{1cm}}p{.45\linewidth}}
    \includegraphics[width=0.45\textwidth]{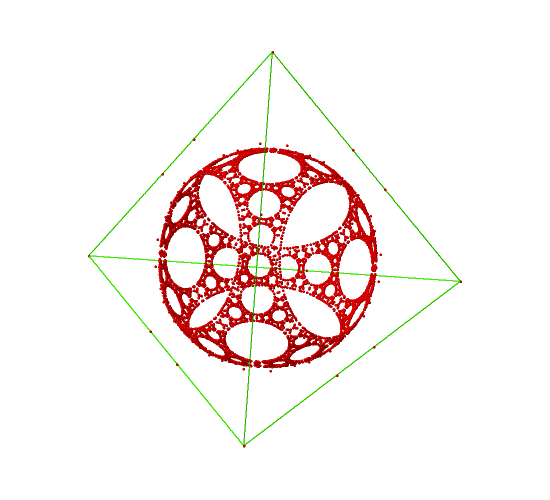}
    &
    \includegraphics[width=0.45\textwidth]{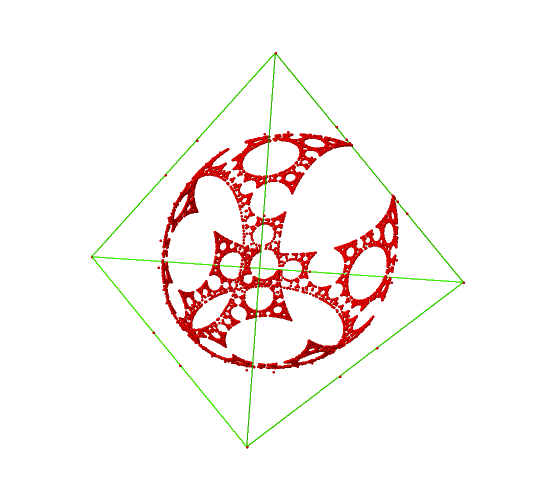}
    \\
    \small (a) Positive roots of depth $\leq$ 7 for the Coxeter system of rank 4 with a complete Coxeter graph with all edges labeled by 4. This Coxeter system is of level~$2$.
    &
    \small (b) Positive roots of depth $\leq$ 7 for the Coxeter system of rank 4 with a complete Coxeter graph with all edges labeled by 4 except one dotted edge labeled by $-1.1$. This Coxeter system is of level~$3$.
  \end{tabular}
  \caption{The pattern of a ball packing and a ball cluster approximated by roots generated by rank-$4$ Coxeter systems, seen in the affine space spanned by simple roots.}
  \label{fig:patterns}
\end{figure}

We say that a Coxeter system is Lorentzian if, in the geometric representation mentioned above, the Coxeter group acts on a Lorentz space as a discrete reflection group generated by reflections in the hyperplanes orthogonal to the basis with respect to the bilinear form; see Section~\ref{ssec:Represent}.  In many examples of Lorentzian Coxeter systems, fractal patterns of ball packings appear while visualizing limit roots on an affine hyperplane; see \cite[Figure~1(b)]{hohlweg_asymptotical_2013}, \cite[Figure~1]{hohlweg_limit_2013} and Figure~\ref{fig:patterns} of the present article. A description of this fractal structure is conjectured in \cite[Section~3.2]{hohlweg_asymptotical_2013} and proved in \cite[Theorem~4.10]{dyer_imaginary2_2013}. In \cite{hohlweg_limit_2013}, Hohlweg, Pr\'eaux and Ripoll prove that the set of limit roots of a Coxeter group~$W$ acting on a Lorentz space is equal to the limit set of $W$ seen as a discrete subgroup of hyperbolic isometries. This explains the pattern of Apollonian disk packing left by the limit roots of the universal Coxeter group of rank~$4$.

While investigating limit roots, we observed that patterns appearing in these examples are similar to the ball packings studied by Boyd and Maxwell, which generalizes the renowed Apollonian ball packings.  One way to generate an Apollonian ball packing is by inversion; see for instance~\cite{graham_apollonian_2005,graham_apollonian_2006}.  In~\cite{boyd_new_1974}, Boyd proposed a class of infinite ball packings generalizing this construction, which is later related to Lorentzian Coxeter systems by Maxwell~\cite{maxwell_sphere_1982}.  Maxwell's approach relies on a correspondance between space-like directions and balls.  More specifically, in the geometric representation of a Coxeter group, weights are vectors ``dual'' to the roots, and the \emph{Boyd--Maxwell ball cluster} refers to the set of balls corresponding to space-like weights.  Maxwell proved that a Boyd--Maxwell ball cluster is a ball packing if and only if the Lorentzian Coxeter system is of ``level~$2$''; see Section~\ref{ssec:bmpacking}.  

The main result of this paper unifies the study of limit roots and the work of Boyd and Maxwell.  Notions involved in Theorem~\ref{thm:main} are formally defined in Section~\ref{sec:Prelim}.  
\begin{theorem}\label{thm:main}
  The set of limit roots of a Lorentzian Coxeter system is the residual set of the corresponding Boyd--Maxwell ball clusters. 
\end{theorem}
Here, the residual set is the complement of the interiors of balls in the cluster.  This theorem implies that~\cite[Theorem~4.10]{dyer_imaginary2_2013} and~\cite[Theorem~1.2]{higashitani_analysis_2013} (see Theorem~\ref{thm:fractal} below) can be deduced from~\cite[Theorem~3.2]{maxwell_sphere_1982} (see Theorem~\ref{thm:maxwell} below) in the Lorentzian case; see Section~\ref{ssec:Main} and~\ref{ssec:cluster}.

We first prove the main result for Lorentzian Coxeter systems of level~$2$. In this case, the Boyd--Maxwell ball cluster is a ball packing, as illustrated in Figure~\ref{fig:patterns}(a).  The proof is based on the study of limit directions of weights, which turn out to coincide with limit roots for a Lorentzian Coxeter system; see Theorem~\ref{thm:limitweight}.  We then extend the same arguments to Lorentzian Coxeter systems of level~$\ge 3$. In this case, balls in the Boyd--Maxwell cluster may overlap, as illustrated in Figure~\ref{fig:patterns}(b).  This completes the proof of Theorem~\ref{thm:main}, since Lorentzian Coxeter systems of level~$1$ have been considered in \cite{dyer_imaginary2_2013, hohlweg_limit_2013}.

For Lorentzian Coxeter systems of level~$2$, we also study the tangency graphs of Boyd--Maxwell ball packings.  In~\cite{chen_apollonian_2013}, the tangency graphs of Apollonian ball packings are compared to the $1$-skeleton of stacked polytopes.  In Theorem~\ref{thm:tangency}, we describe the tangency graph of Boyd--Maxwell ball packing in terms of the corresponding Coxeter complex.  Finally, noticing the importance of Maxwell's work, we use the computer algebra system Sage~\cite{sage} to verify the list of irreducible Coxeter systems of level~$2$, which was manually enumerated by Maxwell in \cite{maxwell_sphere_1982}. We find 326 Coxeter graphs, whereas Maxwell found 323.

This paper is organized as follows. In Section~\ref{sec:Prelim}, we recall the notions of geometric representations of Coxeter system, limit root, Coxeter complex and review the work of Boyd and Maxwell. In Section \ref{sec:Main}, we study the relations between limit roots and Boyd--Maxwell ball clusters through the notion of limit weights, and relate tangency graphs of Boyd--Maxwell ball packings to Coxeter complexes. Finally, in Section~\ref{sec:List}, we describe the algorithm that enumerates all level-$2$ Coxeter graphs. The resulting list is presented in Appendix.

%% file: Prelim.tex
\subsection{Geometric representation of a Coxeter group}\label{ssec:Represent}

Let $(W,S)$ be a finitely generated \Dfn{Coxeter system}, where $S$ is a finite set of generators and the \Dfn{Coxeter group}~$W$ is generated with the relations $(st)^{m_{st}}=e$ where $s,t\in S$, $m_{ss}=1$ and $m_{st}=m_{ts}\geq 2$ or $=\infty$ if $s\neq t$. The cardinality $|S|=n$ is the \Dfn{rank} of the Coxeter system $(W,S)$. For an element $w\in W$, the \Dfn{length}~$\ell(w)$ of $w$ is the smallest natural number $k$ such that $w=s_1s_2\dots s_k$ for $s_i\in S$.  The readers are invited to consult \cite{bourbaki_elements_1968,humphreys_reflection_1992} for more details. We associate a matrix~$B$ to $(W,S)$ as follows:
\[
  B_{st}=
	\begin{cases}
		-\cos(\pi/m_{st}) & \text{if}\quad m_{st}<\infty,\\
		-c_{st} & \text{if}\quad m_{st}=\infty,\\
	\end{cases}
\]
for $s,t\in S$, where $c_{st}$ are chosen arbitrarily with $c_{st}=c_{ts}\geq 1$. We say that the Coxeter system~$(W,S)$ associated with the matrix~$B$ is a \Dfn{geometric Coxeter system}, and denote it by~$(W,S)_B$.

Let $V$ be a real vector space of dimension $n$, equipped with a basis $\Delta=\{\alpha_s\}_{s\in S}$.  The matrix~$B$ defines a bilinear form $\B$ on $V$ by $\B(\alpha_s,\alpha_t)=\alpha_s^\intercal B\alpha_t$ for $s,t\in S$.
For a vector $\alpha\in V$ such that $\B(\alpha,\alpha)\ne 0$, we define the reflection~$\sigma_\alpha$
\begin{equation}\label{eqn:Bref}
  \sigma_\alpha(\x):=\x-2\frac{\B(\x,\alpha)}{\B(\alpha,\alpha)}\alpha,\quad\text{for all }\x\in V.
\end{equation}
The homomorphism $\rho:W\to\mathrm{GL}(V)$ that sends $s\in S$ to $\sigma_{\alpha_s}$ is a faithful \Dfn{geometric representation} of the Coxeter group $W$ as a discrete subgroup of the orthogonal group $O_\B(V)$, i.e., the group of linear transformations of $V$ preserving the bilinear form~$\B$. We refer the readers to \cite[Section~1]{hohlweg_asymptotical_2013} for more details.  In the following, we will write~$w(x)$ in place of $\rho(w)(x)$. 

If the matrix~$B$ is positive definite, we say that $(W,S)_B$ is of \Dfn{finite type}; in this case $W$ is a finite group.  If~$B$ is positive semidefinite but not definite, we say that $(W,S)_B$ is of \Dfn{affine type}. In either case, the group $W$ can be represented as a reflection group in the Euclidean space.  If~$B$ has signature $(n-1,1)$, the pair $(V,\B)$ is an $n$-dimensional \Dfn{Lorentz space}, and we say that $(W,S)_B$ is of \Dfn{Lorentzian type}.  In the present paper, a Coxeter system always comes with an associated matrix $B$. Therefore, we sometimes drop the term ``geometric'', and simply call $(W,S)_B$ a Coxeter system.

The set $Q=\{\x\in V\mid\B(\x,\x)=0\}$ is called \Dfn{isotropic cone}, or \Dfn{light cone} if $(V,\B)$ is a Lorentz space. In a Lorentz space, a vector $\x$ is \Dfn{space-like} (resp.~\Dfn{time-like, light-like}) if $\B(\x,\x)$ is positive (resp.~negative, zero). In \cite{maxwell_sphere_1982}, Maxwell uses the term ``real'' for space-like vectors. The following proposition plays an essential role in the proofs in the present paper; see for instance~\cite[Theorem~2.3]{cecil_lie_2008}.
\begin{proposition}\label{prop:twolights}
  Let $(V,\B)$ be a Lorentz space and $\x,\y\in Q$ be two light-like vectors. Then $\B(\x,\y)=0$ if and only if $\x=c\y$ for some $c\in\mathbb{R}$.
\end{proposition}
Let $\Phi:=W(\Delta)$ be the orbit of $\Delta$ under the action of~$W$. The vectors in $\Delta$ are called \Dfn{simple roots}, and the vectors in $\Phi$ are called \Dfn{roots}. The roots~$\Phi$ are partitioned into \Dfn{positive roots} $\Phi^+=\cone(\Delta)\cap\Phi$ and \Dfn{negative roots} $\Phi^-=-\Phi^+$.  In \cite{hohlweg_asymptotical_2013} and \cite{dyer_imaginary2_2013}, simple roots only need to be positively independent but not necessarily linearly independent.  The \Dfn{depth}~$\depth(\gamma)$ for $\gamma\in\Phi^+$ is the smallest integer $k$ such that $\gamma=s_1s_2\dots s_{k-1}(\alpha)$, for $s_i\in S$ and $\alpha\in\Delta$. 

Let $V^*$ be the dual vector space of $V$ with dual basis $\Delta^*$. If the bilinear form~$\B$ is non-singular, which is the case for Lorentz spaces, $V^*$ can be identified with $V$, and $\Delta^*=\{\omega_s\}_{s\in S}$ can be identified with a set of vectors in $V$ such that
\begin{equation}\label{eqn:dualbasis}
  \B(\alpha_s,\omega_t)=\delta_{st},
\end{equation}
where $\delta_{st}$ is the Kronecker delta function. Vectors in $\Delta^*$ are called \Dfn{fundamental weights}, and vectors in the orbit $\Omega:=W(\Delta^*)$ are called \Dfn{weights}.

\begin{remark}
  In the present article, we are mainly concerned with Coxeter groups acting on Lorentz spaces, therefore we use the term ``Lorentzian''. In the literature, the term \emph{hyperbolic} is used, but with different meanings. In~\cite{bourbaki_elements_1968, humphreys_reflection_1992}, the term \emph{hyperbolic} stands for what we call Lorentzian \emph{of level 1} (see Section~\ref{ssec:bmpacking} for the definition), while \emph{compact} hyperbolic stand for what we call \emph{strict} Lorentzian of level 1 (see Section~\ref{ssec:remarks} for the definition).  In~\cite[Section~9.1]{dyer_imaginary_2013} and~\cite{dyer_imaginary2_2013}, if the simple roots are linearly independent, the term \emph{weakly hyperbolic} corresponds to what we call Lorentzian. Whereas in \cite{vinberg_discrete_1971,maxwell_hyperbolic_1978,maxwell_sphere_1982}, the term \emph{hyperbolic} stands for what we call Lorentzian.
  See~\cite[Section~3.5]{hohlweg_limit_2013} and Remark~3.10 therein for more discussion on terminology.
\end{remark}

\subsection{Limit roots}
As observed in \cite[Section~2.1]{hohlweg_asymptotical_2013}, the set of roots $\Phi\subset V$ is discrete and has no limit point. Nevertheless, it is possible to study the asymptotic directions of the roots. For this, we pass to the \Dfn{projective space} $\mathbb{P}V$, i.e., the space of $1$-dimensional subspaces of~$V$. For a non-zero vector $\x\in V\setminus\{0\}$, let $\proj\x\in\mathbb{P}V$ denote the line passing through $\x$ and the origin. The group action of $W$ on $V$ by reflection induces a \Dfn{projective action} of $W$ on $\mathbb{P}V$:
\[
	w\cdot\proj\x=\proj{w(\x)},\quad w\in W,\quad\x\in V.
\]
One verifies that this is indeed a group action. For a set $X\subset V$, we define the corresponding projective set 
\[\proj X:=\{\proj\x \in \mathbb{P}V\mid\x\in X\}.\]
In this sense, we have \Dfn{projective roots} $\proj\Phi$, \Dfn{projective weights} $\proj\Omega$ and the \Dfn{projective isotropic cone} $\proj Q$.

Let~$h(\x)$ denote the sum of the coordinates of $\x$ in the basis $\Delta$, and call it the \Dfn{height} of the vector~$\x$.  In a Lorentz space, we say that $\x$ is \Dfn{future-directed} (resp. \Dfn{past-directed})\footnotemark if $h(\x)$ is positive (resp. negative). The hyperplane $\{\x\in V\mid h(\x)=1\}$ is the affine subspace $\affine(\Delta)$ spanned by the simple roots. It is useful to identify the projective space $\mathbb{P}V$ with the affine subspace $\affine(\Delta)$ plus a \Dfn{projective hyperplane at infinity}. For a vector $\x\in V$, if $h(\x)\ne 0$, $\proj\x$ is identified with the \Dfn{projective vector}
\[
	\x/h(\x)\in\affine(\Delta).
\]
Otherwise, if $h(\x)=0$, the direction $\proj\x$ is identified with a point on the projective hyperplane at infinity. We avoid the term ``normalized roots'' used in \cite{dyer_imaginary2_2013} and other literature, because the same term is used in \cite{maxwell_sphere_1982} for different objects, which are also important in this paper; see Section \ref{ssec:bmpacking}. For a simple root~$\alpha\in\Delta$, the affine picture of $\proj\alpha$ is $\alpha$ itself. In fact, if $h(\x)\ne 0$, $\proj\x$ is identified with the intersection of $\affine(\Delta)$ with the straight line passing through~$\x$ and the origin. In this sense, the projective roots~$\proj\Phi$, projective weights $\proj\Omega$ and projective isotropic cone $\proj Q$ are respectively identified with the intersection of $\affine(\Delta)$ with the $1$-subspaces spanned by the roots $\Phi$, weights $\Omega$ and isotropic cone~$Q$. In a Lorentz space, the projective light cone~$\proj Q$ is projectively equivalent to a sphere in the affine picture; see for instance~\cite[Proposition~4.13]{dyer_imaginary2_2013}. The affine subspace $\affine(\Delta)$ is practical for visualizing the projective vectors and developing geometric intuitions. In Figure \ref{fig:fundweights}, simple roots, fundamental weights and some positive roots are represented in $\affine(\Delta)$.
\footnotetext{By using the term future- or past-directed, we are assuming that the hyperplane $h(\x)=0$ intersects the light cone only at the origin.  If this is not the case, one can replace $h(\x)$ with any positively weighted sum of the coordinates.  The only requirement is that the hyperplane $h(\x)=1$ is transverse to $\Phi^+$; see~\cite[Section~5.2]{hohlweg_asymptotical_2013}.}

\input{Fund_weights.tex}

\begin{definition}[{Hohlweg--Labb\'e--Ripoll \cite[Definition~2.12]{hohlweg_asymptotical_2013}}]
  The set $E(\Phi)$ of \Dfn{limit roots} is the set of accumulation points of $\proj\Phi$. In other words,
  {\small
  	\[
    E(\Phi)=\{\proj{\x}\in \mathbb{P}V \mid \text{ there is an injective sequence } (\gamma_i)_{i\in\mathbb{N}} \in \Phi \text{ such that } \lim_{i\rightarrow \infty} \proj \gamma_i=\proj{\x}\}.
  \]}
\end{definition}
In~\cite[Theorem~2.7]{hohlweg_asymptotical_2013}, the authors assert that
\[
  E(\Phi)\subseteq\proj Q\cap\cone(\Delta),
\]
see also \cite[Proposition 5.3]{dyer_imaginary_2013}.  Consequently, there is no limit root in the set $\proj Q\setminus\cone(\Delta)$. If the set~$\proj Q\setminus\cone(\Delta)$ consists of open balls (spherical caps), it was conjectured that $E(\Phi)$ is equal to the complement of the $W$-orbit of these balls; see~\cite[Section~3.2]{hohlweg_asymptotical_2013}.  This conjecture is proved in~\cite[Theorem~1.2]{higashitani_analysis_2013} for Lorentzian Coxeter systems, and more generally in \cite[Theorem~4.10]{dyer_imaginary2_2013} for $\Delta$ positively independent.  The present paper relates this result, presented in Theorem~\ref{thm:fractal}, to the result of Maxwell.

\begin{theorem}[{\cite[Theorem~4.10]{dyer_imaginary2_2013}}, {\cite[Theorem~1.2]{higashitani_analysis_2013}}]\label{thm:fractal}
  Let $(W,S)_B$ be an irreducible Lorentzian Coxeter system. Then \[E(\Phi)=\proj Q\setminus(W\cdot(\proj Q\setminus\convex(\Delta))).\] In particular, if $\proj Q\subset\convex(\Delta)$, then $E(\Phi)=\proj Q$.
\end{theorem}

\begin{remark}
  In~\cite{dyer_imaginary2_2013}, the group $W$ acts on the affine space $\affine(\Delta)$, but the action is not defined everywhere.  For this reason, Theorem~4.10 of~\cite{dyer_imaginary2_2013} is stated in terms of $\proj Q_\text{act}$, the part of $\proj Q$ where $W$ acts. This is however not necessary in our setting, because the action of $W$ is well defined on the projective space $\mathbb{P}V$.
\end{remark}

\begin{remark}
	If a Lorentzian Coxeter system is reducible, then one of its irreducible Coxeter subsystems is Lorentzian, and all others are of finite type.  By \cite[Proposition 2.15]{hohlweg_asymptotical_2013}, all the limit roots come from the Lorentzian Coxeter subsystem.  We shall therefore focus on irreducible Lorentzian Coxeter systems.
\end{remark}

The following theorem is useful for the proofs.

\begin{theorem}[{\cite[Theorem 3.1]{dyer_imaginary2_2013}}]\label{thm:minimal}
	The set of limit roots $E_\Phi$ is a minimal set under the action of $W$.  That is, for any limit root $\proj\x\in E_\Phi$, the orbit $W\cdot\proj\x$ is dense in $E_\Phi$.
\end{theorem}

\subsection{Coxeter complex}
For a Lorentzian Coxeter system $(W,S)_B$, let $C=\cone(\Delta^*)$ be the closed cone over fundamental weights. Equivalently, $C$ is the intersection of the half-spaces $\{\x\in V\mid\B(\x,\alpha)\ge 0\}$ for simple roots $\alpha\in\Delta$. The \Dfn{Tits cone}
\[
  T=\cone(\Omega)=\bigcup_{w\in W}w(C)
\]
is the closed cone spanned by weights. It contains one component of the $Q\setminus\{0\}$~\cite[Corollary 1.3]{maxwell_sphere_1982}. By abuse of language, we consider $\proj C$ in the affine picture of $\mathbb{P}V$ as an $(n-1)$-dimensional simplex supported by projective hyperplanes
\[
  \proj H_\alpha=\{\proj\x\in\mathbb{P}V\mid\B(\x,\alpha)=0\}
\]
with $\alpha\in\Delta$. Its vertices are the projective fundamental weights. We call $\proj C$ the \Dfn{fundamental chamber}. The simplices $w\cdot\proj C$ with $w\in W$ are called \Dfn{chambers}. The facets of a chamber are called \Dfn{panels}. For $w\in W$ and $\alpha\in\Delta$, the projective hyperplane $w\cdot\proj H_\alpha$ is called a \Dfn{wall}. The \Dfn{Coxeter complex}~$\C$ associated to the Lorentzian Coxeter system $(W,S)_B$ is the simplicial complex whose maximal simplices correspond to the chambers. The Coxeter complex~$\C$ is a simplicial decomposition of the projective Tits cone $\proj T$ whose vertices correspond to the projective weights. It is pure of dimension $n-1$, where $n$ is the rank of the Coxeter system.

\begin{remark}
  This definition of Coxeter complex is adapted for our purpose.  It applies to Lorentzian Coxeter systems because the Tits cone is strictly convex and does not contain any line through the origin \cite[Section~2.6.3]{abramenko_buildings_2008}.  This is however not true for finite or affine Coxeter groups.  We refer the readers to~\cite[Chapter 3]{abramenko_buildings_2008} for a combinatorial definition in terms of cosets, which applies to general Coxeter systems.
\end{remark}

The group $W$ acts simply transitively on the chambers of $\C$. The \Dfn{dual graph} of $\C$ is the Cayley graph of $(W,S)$. Two chambers are \Dfn{adjacent} if they share a panel. A \Dfn{gallery} is a sequence of chambers $(\proj C_0,\dots,\proj C_k)$ such that consecutive chambers are adjacent, and $k$ is the \Dfn{length} of the gallery.  We say that a gallery~$(\proj C_0,\dots,\proj C_k)$ \Dfn{connects} two simplices $\proj A$ and $\proj A'$ of $\C$ if $\proj A\subseteq\proj C_0$ and $\proj A'\subseteq\proj C_k$.  The \Dfn{gallery distance} $d(\proj A,\proj A')$ between two simplices $\proj A$ and $\proj A'$ is the minimum length of a gallery connecting $\proj A$ and $\proj A'$.  A gallery connecting $\proj A$ and $\proj A'$ with length $d(\proj A,\proj A')$ is called a \Dfn{minimal gallery}.  For an element $w\in W$, its length~$\ell(w)=d(\proj C,w\cdot\proj C)$.  We refer the readers to \cite[Section~1.4.9]{abramenko_buildings_2008} for more details. A pure simplicial complex of dimension $n-1$ is \Dfn{vertex-colorable} if there is a set of $n$ colors and a \Dfn{type function}~$\tau$ that assigns to each vertex of $\C$ a color such that vertices of each chamber have different colors.  The following property~of $\C$ is useful for our purpose~\cite[Theorem~3.5]{abramenko_buildings_2008}:
\begin{theorem}\label{thm:complex}
  The simplicial complex $\C$ is vertex-colorable, and the action of $W$ on $\C$ is type-preserving.
\end{theorem}

In the previous theorem, we can use the fundamental weights $\Delta^*$ as the colors. A vertex $v$ is assigned the color $\omega\in\Delta^*$ if and only if $v$ is in the orbit~$W\cdot\proj\omega$. Correspondingly, every simplex is assigned a type, which is the set of the colors of its vertices. For a panel of type $\Delta^*\setminus\{\omega\}$, we say instead that it is of type $\omega$, to lighten the text.

\subsection{Boyd--Maxwell packing}\label{ssec:bmpacking}

By a \Dfn{ball packing} in a metric space, we mean a set of closed balls with disjoint interiors.  Here, we also regard a closed-half space as a ball of zero curvature, and the complement of an open ball as a ball of negative curvature.  A nice example would be the Apollonian ball packing in dimension $(n-2)$. It is constructed from a set of~$n$ pairwise tangent balls, by repeatedly adding new balls touching $n-1$ pairwise tangent balls. Alternatively, an Apollonian packing can be generated by inversions in the spheres that orthogonally intersect $n-1$ of the initial balls; see~\cite{maxwell_sphere_1982, graham_apollonian_2005, graham_apollonian_2006}.  The group generated by these inversions is called the \Dfn{Apollonian group}.  However, the orbit of the Apollonian group is an infinite ball packing only in dimension $2$ and $3$~\cite{graham_apollonian_2006}. In higher dimensions, restrictions must be imposed to avoid overlap; see~\cite{chen_apollonian_2013}. 

In~\cite{boyd_new_1974}, Boyd presents a new class of infinite ball packings generated by inversions, generalizing Apollonian packings. He characterized these packings in terms of \Dfn{separation} between balls, and explicitly constructed $13$ examples up to dimension nine. Moreover, he noticed a connection to reflection groups. In \cite{maxwell_sphere_1982}, Maxwell revisits these packings, and interprets them using Lorentzian Coxeter groups.

Given a \emph{space-like} vector $\x$ in the Lorentz space $(V,\B)$, the \Dfn{normalized vector}~$\normal{\x}$ of $\x$ is given by
\[
  \normal\x=\x/\sqrt{\B(\x,\x)}
\]
The normalized vector $\normal\x$ lies on the one-sheet hyperboloid $\mathcal{H}=\{\x\in V \mid\B(\x,\x)=1\}$.  Note that $\proj\x=\proj{-\x}$ is the same point in $\mathbb{P}V$, but $\normal\x$ and $\normal{-\x}$ are two different vectors in opposite directions in~$V$.

For $n>2$, there is a classical correspondence between $(n-2)$-dimensional balls and space-like directions in an $n$-dimensional Lorentz space; see for example~\cite[Section~2]{maxwell_sphere_1982}, \cite[Section~2.2]{cecil_lie_2008} or \cite[Section~1.1]{hertrich-jeromin_introduction_2003}. Given a space-like vector~$\x$, let~$H_\x$ be the \Dfn{orthogonal hyperplane} $H_\x=\{\x'\in V\mid\B(\x,\x')=0\}$. In the affine picture of~$\mathbb{P}V$, the intersection of~$\proj Q$ and the half-space $H_\x^-=\{\x'\mid\B(\x,\x')\le 0\}$ is a closed ball (spherical cap) on $\proj Q$. We denote this ball by $\ball(\x)$. After a stereographic projection, $\ball(\x)$ becomes a ball in an $(n-2)$-dimensional Euclidean space. For two space-like vectors $\x$ and $\x'$, if they are \emph{not both} future-directed, we have
\begin{itemize}
  \item $\ball(\x)$ and $\ball(\x')$ are disjoint if $\B(\normal \x,\normal \x')<-1$; 
  \item $\ball(\x)$ is tangent to $\ball(\x')$ if $\B(\normal \x,\normal \x')=-1$;
  \item The boundary of $\ball(\x)$ and $\ball(\x')$ intersect transversally if $\B(\normal \x,\normal \x')>-1$; 
  \item The boundary of $\ball(\x)$ and $\ball(\x')$ intersect transversally at an obtuse angle, or one is contained in the other, if $\B(\normal\x,\normal\x')>0$.
\end{itemize}
In the last case, we say that $\ball(\x)$ and $\ball(\x')$ \Dfn{intersect deeply}.  Therefore, if a set of space-like vectors represents a ball packing, we must have $\B(\normal\x,\normal\x')\le-1$ for any two vectors. The readers are invited to compare with~\cite[Remark~3.2]{hohlweg_limit_2013}.  The packing corresponding to a pair of opposite vectors $\{\x,-\x\}$ is said to be \emph{trivial}; it consists of two balls sharing the same boundary.

\begin{remark}
  One verifies that $-\B(\normal\x,\normal\x')$ is the separation between $\ball(\x)$ and $\ball(\x')$ as defined in \cite{boyd_new_1974}. Given two balls in the Euclidean space of radius~$r$ and $r'$ with centers at distance~$d$ apart, their \Dfn{separation} is defined as $(d^2-r^2-r'^2)/2rr'$.
\end{remark}

To encode geometric Coxeter systems $(W,S)_B$, we adopt Vinberg's convention for Coxeter graphs. That is, if $c_{st}> 1$ the edge~$st$ is dotted and labeled by~$-c_{st}$. This convention is also used by Abramenko--Brown in \cite[Section~10.3.3]{abramenko_buildings_2008} and Maxwell in \cite[Section~1]{maxwell_sphere_1982}.  
A Coxeter graph $G$ is said to be \Dfn{of level~$0$} if it represents a finite or affine Coxeter system. The list of level-$0$ Coxeter graphs can be found in \cite[Chapter~2]{humphreys_reflection_1992}. A graph is \Dfn{of level~$\le r$} if every induced subgraph of $G$ on $n-r$ vertices is of level $0$. A graph is \Dfn{of level~$r$} if it is of level $\le r$ but not of level $\le r-1$. Correspondingly, a Coxeter system~$(W,S)_B$ with a Coxeter graph of level $r$ is said to be \Dfn{of level~$r$}.

For a Lorentzian Coxeter system $(W,S)_B$, while the roots are all space-like, a weight can be space-like, time-like or light-like.  Let $\Omega_r$ be the set of \Dfn{space-like weights}.  We call the set $\{\ball(\omega)\mid\omega\in\Omega_r\}$ the \Dfn{Boyd--Maxwell ball cluster} generated by $(W,S)_B$. Maxwell proved that Coxeter systems of level $2$ are Lorentzian \cite[Proposition~1.6]{maxwell_sphere_1982} and the following theorem. 
\begin{theorem}[{Maxwell \cite[Theorem~3.2]{maxwell_sphere_1982}}]\label{thm:maxwell} 
  Let $(W,S)_B$ be a Lorentzian Coxeter system. The Boyd--Maxwell ball cluster generated by $(W,S)_B$ is a ball packing if and only if $(W,S)_B$ is of level $2$.
\end{theorem}
For example, the Apollonian circle packing is the Boyd--Maxwell ball packing generated by the universal Coxeter system of rank $4$.  Maxwell manually enumerated the Coxeter graphs representing irreducible Coxeter systems of level~$2$, and suggested a computer verification.

\begin{remark}
The Boyd--Maxwell ball packing is trivial if the level-$2$ Coxeter system is reducible.  More generally, the Boyd--Maxwell ball cluster covers the projective light cone if the Lorentzian Coxeter system is reducible.  This gives another reason for focusing on irreducible Coxeter systems.
\end{remark}

\begin{figure}[!ht]
  \begin{center}
    \begin{tikzpicture}[sommet/.style={inner sep=2pt,circle,draw=black,fill=blue,thick,anchor=base}]

      \node[inner sep=0pt] at (0,0) {\includegraphics[width=10cm]{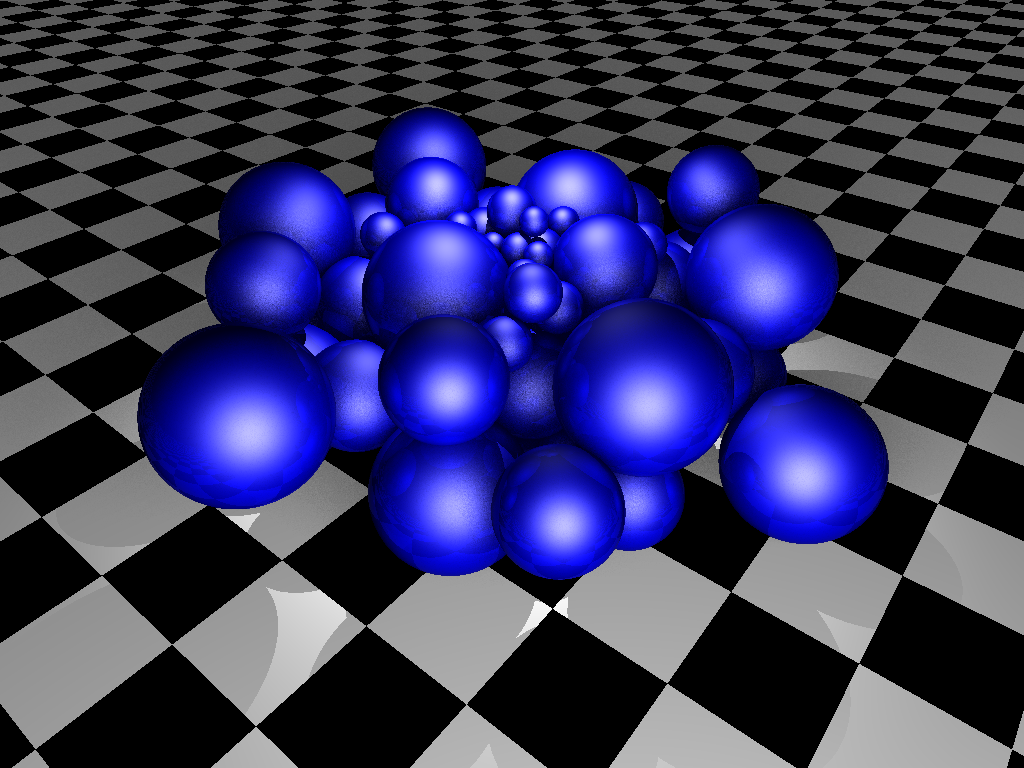}};

      \coordinate (ancre) at (-4.25,-3.1);
      \fill[white] ($(ancre)+(-1,-1)$) rectangle  ($(ancre)+(1,1)$);
      \node[sommet] (alpha) at ($(ancre)+(0:0.6)$) {};
      \node[sommet] (beta) at ($(ancre)+(72:0.6)$) {} edge[thick] node[auto] {4} (alpha);
      \node[sommet] (delta) at ($(ancre)+(144:0.6)$) {} edge[thick] node[auto] {4} (beta);
      \node[sommet] (gamma) at ($(ancre)+(216:0.6)$) {} edge[thick] node[auto] {4} (delta);
      \node[sommet] (eta) at ($(ancre)+(288:0.6)$) {} edge[thick] node[auto] {4} (gamma) edge[thick] node[auto,swap] {4} (alpha);

    \end{tikzpicture}
    \caption{Some balls in the ball packing generated by the level-$2$ Lorentzian Coxeter group whose Coxeter graph is a 5-cycle with all edges labeled by $4$. The coordinates are calculated by Sage \cite{sage} with the help of CHEVIE package~\cite{chevie_1996}, the image is rendered by POV-Ray.}
    \label{fig:ballpack}
  \end{center}
\end{figure}

The reflection $\sigma_\alpha$ with respect to a root $\alpha$ correspond to the inversion with respect to the boundary of $\ball(\alpha)$. These inversions induce a representation of the Lorentzian Coxeter group~$W$ as a subgroup of M\"obius transformations. A Boyd--Maxwell ball packing is therefore generated by inversions from the balls corresponding to space-like fundamental weights. Figure \ref{fig:ballpack} shows an image of a ball packing generated in this way.

The \Dfn{residual set} of a ball packing is the complement of the interiors of all balls in the packing.  The Hausdorff dimension of the residual set of Apollonian disk packings was studied in~\cite{boyd_residual_1973} and calculated in~\cite{mcmullen_hausdorff_1998}; see also~\cite{graham_apollonian_2003}.  The notion of residual set naturally extends to any collection of balls, not necessarily a packing.

%% file: Fund_weights.tex
\begin{figure}[!ht]
  \centering
  \scalebox{1}{
    \begin{tikzpicture}
      [scale=1.5,
	q/.style={black,line join=round,thin},
	racine/.style={red},
	poid/.style={blue},
	racinesimple/.style={black},
	racinedih/.style={blue},
	sommet/.style={inner sep=2pt,circle,draw=black,fill=blue,thick,anchor=base},
      rotate=0]

      \def\grosseur{0.0125}
      \def\grosseursimple{0.025}

      \def\grosseurdih{0.0075}


      \draw[q] (3.32,1.34) -- 
      (3.3,1.44) -- 
      (3.29,1.48) -- 
      (3.28,1.52) -- 
      (3.27,1.55) -- 
      (3.26,1.58) -- 
      (3.24,1.63) -- 
      (3.22,1.68) -- 
      (3.21,1.7) -- 
      (3.2,1.72) -- 
      (3.19,1.74) -- 
      (3.18,1.76) -- 
      (3.17,1.78) -- 
      (3.14,1.83) -- 
      (3.08,1.92) -- 
      (3.05,1.96) -- 
      (3.01,2.01) -- 
      (2.84,2.18) -- 
      (2.83,2.19) -- 
      (2.78,2.23) -- 
      (2.71,2.28) -- 
      (2.68,2.3) -- 
      (2.58,2.36) -- 
      (2.56,2.37) -- 
      (2.54,2.38) -- 
      (2.5,2.4) -- 
      (2.48,2.41) -- 
      (2.41,2.44) -- 
      (2.33,2.47) -- 
      (2.3,2.48) -- 
      (2.23,2.5) -- 
      (2.19,2.51) -- 
      (2.15,2.52) -- 
      (2.1,2.53) -- 
      (2.03,2.54) -- 
      (1.94,2.55) -- 
      (1.93,2.55) -- 
      (1.72,2.55) -- 
      (1.71,2.55) -- 
      (1.62,2.54) -- 
      (1.56,2.53) -- 
      (1.39,2.49) -- 
      (1.36,2.48) -- 
      (1.33,2.47) -- 
      (1.28,2.45) -- 
      (1.23,2.43) -- 
      (1.21,2.42) -- 
      (1.19,2.41) -- 
      (1.17,2.4) -- 
      (1.15,2.39) -- 
      (1.13,2.38) -- 
      (1.08,2.35) -- 
      (0.99,2.29) -- 
      (0.95,2.26) -- 
      (0.88,2.2) -- 
      (0.87,2.19) -- 
      (0.86,2.18) -- 
      (0.85,2.17) -- 
      (0.84,2.16) -- 
      (0.83,2.15) -- 
      (0.82,2.14) -- 
      (0.81,2.13) -- 
      (0.75,2.06) -- 
      (0.72,2.02) -- 
      (0.7,1.99) -- 
      (0.68,1.96) -- 
      (0.65,1.91) -- 
      (0.59,1.79) -- 
      (0.57,1.74) -- 
      (0.56,1.71) -- 
      (0.55,1.68) -- 
      (0.54,1.65) -- 
      (0.53,1.61) -- 
      (0.51,1.51) -- 
      (0.5,1.43) -- 
      (0.5,1.42) -- 
      (0.5,1.24) -- 
      (0.5,1.23) -- 
      (0.51,1.15) -- 
      (0.52,1.09) -- 
      (0.53,1.04) -- 
      (0.55,0.97) -- 
      (0.56,0.94) -- 
      (0.57,0.91) -- 
      (0.58,0.88) -- 
      (0.61,0.81) -- 
      (0.66,0.71) -- 
      (0.72,0.61) -- 
      (0.77,0.54) -- 
      (0.81,0.49) -- 
      (0.82,0.48) -- 
      (0.99,0.31) -- 
      (1.0,0.3) -- 
      (1.05,0.26) -- 
      (1.12,0.21) -- 
      (1.15,0.19) -- 
      (1.2,0.16) -- 
      (1.25,0.13) -- 
      (1.27,0.12) -- 
      (1.29,0.11) -- 
      (1.31,0.1) -- 
      (1.33,0.09) -- 
      (1.35,0.08) -- 
      (1.42,0.05) -- 
      (1.6,-0.01) -- 
      (1.64,-0.02) -- 
      (1.69,-0.03) -- 
      (1.74,-0.04) -- 
      (1.8,-0.05) -- 
      (1.81,-0.05) -- 
      (1.9,-0.06) -- 
      (1.92,-0.06) -- 
      (2.09,-0.06) -- 
      (2.1,-0.06) -- 
      (2.11,-0.06) -- 
      (2.2,-0.05) -- 
      (2.27,-0.04) -- 
      (2.32,-0.03) -- 
      (2.36,-0.02) -- 
      (2.4,-0.01) -- 
      (2.52,0.03) -- 
      (2.57,0.05) -- 
      (2.73,0.13) -- 
      (2.78,0.16) -- 
      (2.81,0.18) -- 
      (2.85,0.21) -- 
      (2.89,0.24) -- 
      (2.9,0.25) -- 
      (3.06,0.41) -- 
      (3.07,0.42) -- 
      (3.1,0.46) -- 
      (3.16,0.55) -- 
      (3.19,0.6) -- 
      (3.2,0.62) -- 
      (3.21,0.64) -- 
      (3.22,0.66) -- 
      (3.23,0.68) -- 
      (3.25,0.73) -- 
      (3.29,0.85) -- 
      (3.3,0.89) -- 
      (3.31,0.93) -- 
      (3.32,0.99) -- 
      (3.33,1.08) -- 
      (3.33,1.09) -- 
      (3.33,1.23) -- 
      (3.33,1.24) -- 
      (3.33,1.25) -- 
      cycle;

      \node[label=left :{$\alpha_s$}] (a) at (0.000000000000000,0.000000000000000) {};
      \fill[racinesimple] (0.000000000000000,0.000000000000000) circle (\grosseursimple);\node[label=right :{$\alpha_t$}] (b) at (4.00000000000000,0.000000000000000) {};
      \fill[racinesimple] (4.00000000000000,0.000000000000000) circle (\grosseursimple);\node[label=above :{$\alpha_r$}] (g) at (2.00000000000000,3.46410161513775) {};
      \fill[racinesimple] (2.00000000000000,3.46410161513775) circle (\grosseursimple);
      \draw[green!75!black] (a) -- (b) -- (g) -- (a);
      \fill[racine] (0.500000000000000,0.866025403784439) circle (\grosseursimple);
      \fill[racine] (0.545454545454545,0.944754985946660) circle (\grosseursimple);
      \fill[racine] (0.551724137931034,0.955614238658691) circle (\grosseursimple);
      \fill[racine] (0.552631578947368,0.957185972603853) circle (\grosseursimple);
      \fill[racine] (0.552763819095477,0.957415019259178) circle (\grosseursimple);
      \fill[racine] (0.557921328040419,0.941345548971032) circle (\grosseur);
      \fill[racine] (0.562779281263767,0.926798160107659) circle (\grosseur);
      \fill[racine] (0.565969062784349,0.917246196546940) circle (\grosseur);
      \fill[racine] (0.567650050864700,0.912718533388279) circle (\grosseur);
      \fill[racine] (0.582294960014878,0.873273152188809) circle (\grosseur);
      \fill[racine] (0.583673469387755,0.869560201350906) circle (\grosseur);
      \fill[racine] (0.586489812662382,0.862670093181326) circle (\grosseur);
      \fill[racine] (0.593582887700534,0.842869644324961) circle (\grosseur);
      \fill[racine] (0.594695843491705,0.842594259185407) circle (\grosseur);
      \fill[racine] (0.601010101010101,0.827146597778403) circle (\grosseur);
      \fill[racine] (0.603336236350238,0.822222283274138) circle (\grosseur);
      \fill[racine] (0.635632360865214,0.753852970604959) circle (\grosseur);
      \fill[racine] (0.637752376882812,0.749365000603990) circle (\grosseur);
      \fill[racine] (0.644355644355644,0.735386207009763) circle (\grosseur);
      \fill[racine] (0.646072426711066,0.733832550629021) circle (\grosseur);
      \fill[racine] (0.653379139043262,0.720191869848972) circle (\grosseur);
      \fill[racine] (0.656518240343348,0.714331576351167) circle (\grosseur);
      \fill[racine] (0.657713836786801,0.712261883206783) circle (\grosseur);
      \fill[racine] (0.681818181818182,0.656079851351847) circle (\grosseur);
      \fill[racine] (0.685561740315801,0.664054467288988) circle (\grosseur);
      \fill[racine] (0.687185138539043,0.661244207864378) circle (\grosseur);
      \fill[racine] (0.691612453805920,0.654143301499574) circle (\grosseur);
      \fill[racine] (0.700000000000000,0.640690638527905) circle (\grosseur);
      \fill[racine] (0.706896551724138,0.627121844119766) circle (\grosseur);
      \fill[racine] (0.712478920741990,0.620675879609421) circle (\grosseur);
      \fill[racine] (0.713836477987421,0.619108307527240) circle (\grosseur);
      \fill[racine] (0.962059620596206,0.332484911118687) circle (\grosseur);
      \fill[racine] (0.963636363636364,0.330664245081331) circle (\grosseur);
      \fill[racine] (0.970149253731343,0.323143807382253) circle (\grosseur);
      \fill[racine] (0.981203007518797,0.315806256267258) circle (\grosseur);
      \fill[racine] (1.00000000000000,0.288675134594813) circle (\grosseur);
      \fill[racine] (1.00279329608939,0.297545040965045) circle (\grosseur);
      \fill[racine] (1.00990099009901,0.291533304244266) circle (\grosseur);
      \fill[racine] (1.01535312180143,0.287197779760653) circle (\grosseur);
      \fill[racine] (1.04010279790339,0.267516885652054) circle (\grosseur);
      \fill[racine] (1.04419556032459,0.264262328966481) circle (\grosseur);
      \fill[racine] (1.05204442935442,0.258333231472780) circle (\grosseur);
      \fill[racine] (1.06419369828072,0.249155578021367) circle (\grosseur);
      \fill[racine] (1.06565176022835,0.247200400699649) circle (\grosseur);
      \fill[racine] (1.07255520504732,0.242839229943060) circle (\grosseur);
      \fill[racine] (1.07352157698455,0.242228735741519) circle (\grosseur);
      \fill[racine] (1.25000000000000,0.01) circle (\grosseursimple);
      \fill[racine] (1.26009693053312,0.124361931974071) circle (\grosseur);
      \fill[racine] (1.26185958254269,0.123248397123023) circle (\grosseur);
      \fill[racine] (1.27490039840637,0.115010013782794) circle (\grosseur);
      \fill[racine] (1.27519454247449,0.116368607038058) circle (\grosseur);
      \fill[racine] (1.29105835393293,0.108184128543955) circle (\grosseur);
      \fill[racine] (1.30416451960463,0.101422378294894) circle (\grosseur);
      \fill[racine] (1.30897191164440,0.0991902562283786) circle (\grosseur);
      \fill[racine] (1.37718108247182,0.0675200296100512) circle (\grosseur);
      \fill[racine] (1.38778460426455,0.0625967042851058) circle (\grosseur);
      \fill[racine] (1.40191553551118,0.0570461628107106) circle (\grosseur);
      \fill[racine] (1.44723618090452,2.50668659587858) circle (\grosseursimple);
      \fill[racine] (1.44736842105263,2.50691564253390) circle (\grosseursimple);
      \fill[racine] (1.44827586206897,2.50848737647906) circle (\grosseursimple);
      \fill[racine] (1.45454545454545,2.51934662919109) circle (\grosseursimple);
      \fill[racine] (1.45695364238411,0.01) circle (\grosseursimple);
      \fill[racine] (1.46218526254133,0.0333725907231712) circle (\grosseur);
      \fill[racine] (1.46321525885559,2.51076574830148) circle (\grosseur);
      \fill[racine] (1.48208955223881,2.51535139666346) circle (\grosseur);
      \fill[racine] (1.49227799227799,2.51782675308758) circle (\grosseur);
      \fill[racine] (1.49869294188619,2.51918831513084) circle (\grosseur);
      \fill[racine] (1.50000000000000,2.59807621135332) circle (\grosseursimple);
      \fill[racine] (1.52260111022998,0.01) circle (\grosseursimple);
      \fill[racine] (1.52908157612910,2.52563825058771) circle (\grosseur);
      \fill[racine] (1.53431708991078,2.52674947940127) circle (\grosseur);
      \fill[racine] (1.54279063398603,2.52830877159433) circle (\grosseur);
      \fill[racine] (1.54714738510301,0.01) circle (\grosseursimple);
      \fill[racine] (1.55687719254127,0.01) circle (\grosseursimple);
      \fill[racine] (1.56281407035176,2.53279791458564) circle (\grosseur);
      \fill[racine] (1.56281407035176,2.53279791458565) circle (\grosseur);
      \fill[racine] (1.56347170601476,2.53211447955693) circle (\grosseur);
      \fill[racine] (1.57627118644068,2.53446982576463) circle (\grosseur);
      \fill[racine] (1.57922292642601,2.53483654973775) circle (\grosseur);
      \fill[racine] (1.69751243781095,2.54953283051100) circle (\grosseur);
      \fill[racine] (1.70024436600597,2.54987224508239) circle (\grosseur);
      \fill[racine] (1.71257485029940,2.55140418360445) circle (\grosseur);
      \fill[racine] (1.71499821873887,2.55087211916271) circle (\grosseur);
      \fill[racine] (1.72746721101834,2.55171714740135) circle (\grosseur);
      \fill[racine] (1.73514934933585,2.55223776876881) circle (\grosseur);
      \fill[racine] (1.73748697852076,2.55232408798259) circle (\grosseur);
      \fill[racine] (1.79905078616860,2.55459739099737) circle (\grosseur);
      \fill[racine] (1.80311908410709,2.55474761682520) circle (\grosseur);
      \fill[racine] (1.80487804878049,2.56287192664663) circle (\grosseur);
      \fill[racine] (1.81182998318327,2.55475533088899) circle (\grosseur);
      \fill[racine] (1.83633640886914,2.55477703291184) circle (\grosseur);
      \fill[racine] (1.85776235906331,2.55676537249109) circle (\grosseur);
      \fill[racine] (1.87468776019983,2.55481099551063) circle (\grosseur);
      \fill[racine] (1.88103508840404,2.55407806918116) circle (\grosseur);
      \fill[racine] (2.12261152902369,2.52618322457684) circle (\grosseur);
      \fill[racine] (2.12936499409047,2.52540340160198) circle (\grosseur);
      \fill[racine] (2.14748201438849,2.52331142829279) circle (\grosseur);
      \fill[racine] (2.16583703731904,2.51679591871416) circle (\grosseur);
      \fill[racine] (2.19849779600575,2.50908790740979) circle (\grosseur);
      \fill[racine] (2.20512820512820,2.51665501954452) circle (\grosseur);
      \fill[racine] (2.20721178263078,2.50703138728862) circle (\grosseur);
      \fill[racine] (2.21314514147888,2.50541418152704) circle (\grosseur);
      \fill[racine] (2.25751837825796,2.49331974125865) circle (\grosseur);
      \fill[racine] (2.26112163258996,2.49233763250451) circle (\grosseur);
      \fill[racine] (2.27030810580566,2.48956983131642) circle (\grosseur);
      \fill[racine] (2.28211730303932,2.48601182713284) circle (\grosseur);
      \fill[racine] (2.28471411901984,2.48590722673246) circle (\grosseur);
      \fill[racine] (2.29268292682927,2.48282850036617) circle (\grosseur);
      \fill[racine] (2.29379265518362,2.48239975995409) circle (\grosseur);
      \fill[racine] (2.44312280745873,0.01) circle (\grosseursimple);
      \fill[racine] (2.45285261489699,0.01) circle (\grosseursimple);
      \fill[racine] (2.47739888977002,0.01) circle (\grosseursimple);
      \fill[racine] (2.48449385605617,2.40872289018258) circle (\grosseur);
      \fill[racine] (2.48609680741504,2.40810359445724) circle (\grosseur);
      \fill[racine] (2.49708656840976,2.40282267076862) circle (\grosseur);
      \fill[racine] (2.49785407725322,2.40356120649758) circle (\grosseur);
      \fill[racine] (2.51411420538733,2.39464035893133) circle (\grosseur);
      \fill[racine] (2.52584487683933,2.38900340460775) circle (\grosseur);
      \fill[racine] (2.53162895472794,2.38599514970746) circle (\grosseur);
      \fill[racine] (2.53755260602137,0.0350447314577711) circle (\grosseur);
      \fill[racine] (2.54304635761589,0.01) circle (\grosseursimple);
      \fill[racine] (2.57018498367791,2.36594245241363) circle (\grosseur);
      \fill[racine] (2.57142857142857,2.47435829652697) circle (\grosseursimple);
      \fill[racine] (2.57960644007156,2.36104242462877) circle (\grosseur);
      \fill[racine] (2.58180706977111,0.0533550733755174) circle (\grosseur);
      \fill[racine] (2.59094437257439,2.35451356868483) circle (\grosseur);
      \fill[racine] (2.59605911330049,0.0592518749168335) circle (\grosseur);
      \fill[racine] (2.60330058558777,0.0628678365037341) circle (\grosseur);
      \fill[racine] (2.62919652961147,2.33248637609238) circle (\grosseur);
      \fill[racine] (2.64516129032258,2.34664948122235) circle (\grosseursimple);
      \fill[racine] (2.66141732283465,2.31849320698196) circle (\grosseursimple);
      \fill[racine] (2.66536203522505,2.31166076470054) circle (\grosseursimple);
      \fill[racine] (2.66634098680997,2.30996517081731) circle (\grosseursimple);
      \fill[racine] (2.69430635024857,0.108310715540686) circle (\grosseur);
      \fill[racine] (2.69719499831024,0.109753134984510) circle (\grosseur);
      \fill[racine] (2.70738223933694,0.115565893098722) circle (\grosseur);
      \fill[racine] (2.72304439746300,0.124502595047234) circle (\grosseur);
      \fill[racine] (2.72340425531915,0.122840482806303) circle (\grosseur);
      \fill[racine] (2.74030406453615,0.134350822802426) circle (\grosseur);
      \fill[racine] (2.74388867995487,0.136792278897881) circle (\grosseur);
      \fill[racine] (2.75000000000000,0.01) circle (\grosseursimple);
      \fill[racine] (2.85925774221390,0.215369327584304) circle (\grosseur);
      \fill[racine] (2.86140089418778,0.216829013168086) circle (\grosseur);
      \fill[racine] (2.87084300077339,0.223259964883846) circle (\grosseur);
      \fill[racine] (2.87187954849449,0.225175008422413) circle (\grosseur);
      \fill[racine] (2.88756752173748,0.237670100483247) circle (\grosseur);
      \fill[racine] (2.89359179169617,0.242468286110368) circle (\grosseur);
      \fill[racine] (2.89745574794985,0.245730882051304) circle (\grosseur);
      \fill[racine] (2.91820933908431,0.263254523042501) circle (\grosseur);
      \fill[racine] (2.92225700025379,0.266672233294591) circle (\grosseur);
      \fill[racine] (2.92891419796532,0.272624287326089) circle (\grosseur);
      \fill[racine] (2.93129770992366,0.264435237796775) circle (\grosseur);
      \fill[racine] (2.94051261896073,0.282994180400631) circle (\grosseur);
      \fill[racine] (2.94736842105263,0.287875480482362) circle (\grosseur);
      \fill[racine] (2.94957983193277,0.291100976061996) circle (\grosseur);
      \fill[racine] (2.94989979959920,0.291567671013599) circle (\grosseur);
      \fill[racine] (3.25298421255295,0.733637230776190) circle (\grosseur);
      \fill[racine] (3.25327951564077,0.734067950735548) circle (\grosseur);
      \fill[racine] (3.25531914893617,0.737042896837820) circle (\grosseur);
      \fill[racine] (3.25752336315052,0.745468410790983) circle (\grosseur);
      \fill[racine] (3.26167220712227,0.756613658085414) circle (\grosseur);
      \fill[racine] (3.26453863067708,0.764313875202886) circle (\grosseur);
      \fill[racine] (3.26572528883184,0.767823978023687) circle (\grosseur);
      \fill[racine] (3.27007299270073,0.758562397475421) circle (\grosseur);
      \fill[racine] (3.27711254683406,0.801507180110648) circle (\grosseur);
      \fill[racine] (3.27829208388438,0.804996218896550) circle (\grosseur);
      \fill[racine] (3.28016235456987,0.811177913876863) circle (\grosseur);
      \fill[racine] (3.28596333420256,0.830351548449951) circle (\grosseur);
      \fill[racine] (3.28741092636580,0.831969509046017) circle (\grosseur);
      \fill[racine] (3.29082086766775,0.846406866255772) circle (\grosseur);
      \fill[racine] (3.29194088642621,0.851148915888199) circle (\grosseur);
      \fill[racine] (3.31069979214095,0.930572283824989) circle (\grosseur);
      \fill[racine] (3.31218697829716,0.936868884227573) circle (\grosseur);
      \fill[racine] (3.31555986427533,0.957119690076840) circle (\grosseur);
      \fill[racine] (3.31707317073171,0.957556544021818) circle (\grosseur);
      \fill[racine] (3.31877200287694,0.976405368751778) circle (\grosseur);
      \fill[racine] (3.32023121387283,0.985166470894668) circle (\grosseur);
      \fill[racine] (3.32056194125160,0.988057931946486) circle (\grosseur);
      \fill[racine] (3.32953249714937,1.06648510386225) circle (\grosseur);
      \fill[racine] (3.33014354066986,1.07182729399797) circle (\grosseur);
      \fill[racine] (3.33084592329806,1.08442547350144) circle (\grosseur);
      \fill[racine] (3.33226810192583,1.10993416480320) circle (\grosseur);
      \fill[racine] (3.33365901319003,1.15413644432044) circle (\grosseursimple);
      \fill[racine] (3.33463796477495,1.15244085043722) circle (\grosseursimple);
      \fill[racine] (3.33858267716535,1.14560840815579) circle (\grosseursimple);
      \fill[racine] (3.35483870967742,1.11745213391540) circle (\grosseursimple);
      \fill[racine] (3.42857142857143,0.989743318610787) circle (\grosseursimple);

      \node[label=right :{$\omega_s$}] at (3.33806146572104, 1.88355406969665) {};
      \node[label=left :{$\omega_t$}] at (0.172972972972973, 2.17208533705935) {};
      \node[label=below :{$\omega_r$}] at (2.00898472596586, -0.130720815665576) {};

      \fill[poid] (3.33806146572104, 1.88355406969665) circle (\grosseursimple);
      \fill[poid] (0.172972972972973, 2.17208533705935) circle (\grosseursimple);
      \fill[poid] (2.00898472596586, -0.130720815665576) circle (\grosseursimple);

      \draw[thin,blue] (3.33806146572104, 1.88355406969665) -- (3.33333333333333, 1.15470053837925);
      \draw[thin,blue] (3.33806146572104, 1.88355406969665) -- (2.66666666666667, 2.30940107675850);

      \draw[thin,blue] (0.172972972972973, 2.17208533705935) -- (0.552786404500042, 0.957454138327394);
      \draw[thin,blue] (0.172972972972973, 2.17208533705935) -- (1.44721359549996, 2.50664747681036);

      \draw[thin,blue] (2.00898472596586, -0.130720815665576) -- (2.43643578047198, 0.000000000000000);
      \draw[thin,blue] (2.00898472596586, -0.130720815665576) -- (1.56356421952802, 0.000000000000000);

      \coordinate (ancre) at (-1.5,2.6);
      \node[sommet,label=below left:$s$] (alpha) at (ancre) {};
      \node[sommet,label=below right :$t$] (beta) at ($(ancre)+(0.5,0)$) {} edge[thick,dotted] node[auto] {$-1.1$} (alpha);
      \node[sommet,label=above:$r$] (gamma) at ($(ancre)+(0.25,0.43)$) {} edge[thick,dotted] node[auto,swap] {$-1.5$} (alpha) edge[thick,dotted] node[auto] {$-1.25$} (beta);

    \end{tikzpicture}
  }
  \caption{\label{fig:fundweights} Simple roots, fundamental weights, and positive roots of depth $\leq 6$ of a geometric Coxeter system of rank~$3$ seen in the affine space spanned by the simple roots. The Coxeter graph is shown in the upper-left corner.}

\end{figure}

%% file: Main.tex
\subsection{Limit weights}
Let $(W,S)_B$ be a (not necessarily Lorentzian) geometric Coxeter system as described in Section \ref{ssec:Represent}. When the bilinear form $\B$ is non-singular, we define the set of limit weights $E(\Omega)$ analogously to limit roots.
\begin{definition}
  The set of \Dfn{limit weights} $E(\Omega)$ is the set of accumulation points of the projective weights $\proj\Omega$. That is
  {\small\[
    E(\Omega)=\{\proj{\x}\in \mathbb{P}V \mid \text{there is an injective sequence } (\omega_i)_{i\in\mathbb{N}} \in \Omega\text{ such that } \lim_{i\rightarrow \infty} \proj \omega_i=\proj{\x}\}.
  \]}
\end{definition}

We recall the following theorem about limit roots.

\begin{theorem}[{Hohlweg--Labb\'e--Ripoll \cite[Theorem~2.7]{hohlweg_asymptotical_2013}}]\label{thm:ConvRoot}
  Consider an injective sequence of roots $(\gamma_k)_{k\in\mathbb{N}}$ and suppose that $(\proj\gamma_k)_{k\in\mathbb{N}}$ converges to a limit $\proj\beta$. Then
  \begin{enumerate}
      \renewcommand{\labelenumi}{(\roman{enumi})}
      \renewcommand{\theenumi}{\ref{thm:ConvRoot}\labelenumi}
    \item \label{thm:ConvRoot1} $h(\gamma_k)$ tends to $+\infty$,
    \item \label{thm:ConvRoot2} $\proj\beta$ lies in $\proj Q$.
  \end{enumerate}
\end{theorem}

\begin{remark}\label{rem:growthdepth}
	Theorem~\ref{thm:ConvRoot1} is not in the statement of \cite[Theorem~2.7]{hohlweg_asymptotical_2013}, but mentioned in its proof.  It is proved in \cite[Lemma~2.10]{hohlweg_asymptotical_2013} that the squared Euclidean norm of a positive root grows at least linearly with its depth.  Then, since the height of a positive root is nothing but its $L_1$-norm, Theorem~\ref{thm:ConvRoot1} follows from the equivalence of the norms.
\end{remark}

Here is an analogous result for limit weights.
\begin{theorem}\label{thm:ConvWeight}
  Consider an injective sequence of weights $(\omega_k)_{k\in\mathbb{N}}$ and suppose that $(\proj\omega_k)_{k\in\mathbb{N}}$ converges to a limit $\proj\psi$. Then
  \begin{enumerate}
      \renewcommand{\labelenumi}{(\roman{enumi})}
      \renewcommand{\theenumi}{\ref{thm:ConvWeight}\labelenumi}
    \item \label{thm:ConvWeight1} $h(\omega_k)$ tends to $-\infty$,
    \item \label{thm:ConvWeight2} $\proj\psi$ lies in $\proj Q$.
  \end{enumerate}
\end{theorem}

In preparation for the proof, we make the following observations.  For $\alpha_t\in\Delta$ and $\omega_s\in\Delta^*$, we have from Equation \eqref{eqn:dualbasis} that
\begin{equation}\label{eqn:WeightSimpleRef}
  \sigma_{\alpha_t}(\omega_s)=\begin{cases}
    \omega_s & \text{if}\;s\ne t,\\
    \omega_s-2\alpha_s & \text{if}\;s=t.
  \end{cases}
\end{equation}
Let $s\in S$ and $w=s_1s_2\dots s_k\in W$ where $k=\ell(w)$. Define the set
\[
  \prefix^{s}(w):=\{s_1s_2\dots s_m\mid s_{m+1}=s\}.
\]
For any $w'\in\prefix^{s}(w)$, since the expression $w=s_1s_2\dots s_k$ is reduced, we have $\ell(w')<\ell(w's)$ and $w'(\alpha_s)\in\Phi^+$; see for example~\cite[Theorem~5.4]{humphreys_reflection_1992}. From Equation~\eqref{eqn:WeightSimpleRef} we have
\begin{equation}\label{eqn:WeightRef}
  w(\omega_s)=\omega_s-2\cdot\sum_{w'\in\prefix^{s}(w)}w'(\alpha_s).
\end{equation}

\begin{proof}[Proof of Theorem~\ref{thm:ConvWeight}]
  Every weight in the sequence can be written in the form of an element of~$W$
  acting on a fundamental weight. By passing to a subsequence if necessary, we
  may assume that $\omega_k=w_k(\omega_s)$ for a fixed fundamental weight
  $\omega_s\in\Delta^*$ and an injective sequence of elements
  $(w_k)_{k\in\mathbb{N}}$ with increasing length.  

  (i) Using the linearity of $h$ in Equation \eqref{eqn:WeightRef}, we get
	\begin{equation}\label{eqn:heightgrowth}
    h(\omega_k)=h(w_k(\omega_s))=h(\omega_s)-2\cdot\sum_{w'\in\prefix^{s}(w_k)}h(w'(\alpha_s))
  \end{equation}
	While $h(\omega_s)$ remains constant, we claim that the summation in~\eqref{eqn:heightgrowth} diverges as~$k$ tends to~$+\infty$.  To prove this claim, we first notice that all the summands are positive.  Since $S$ is finite, there are only finitely many positive roots with a bounded depth.  As the height of a positive root grows with depth (see Remark~\ref{rem:growthdepth}), there are only finitely many positive roots with a bounded height.  If the summation in~\eqref{eqn:heightgrowth} is bounded by a positive number for infinitely many $k\in\mathbb{N}$, the summation in~\eqref{eqn:WeightRef} contains a bounded number of items chosen from finitely many positive roots for these $k$.  Consequently, the sequence $(\omega_k)$ is not injective as assumed.  This contradiction proves our claim.

  (ii) Since $w_k$ preserves the bilinear form, $\B(\omega_k,\omega_k)=\B(\omega_s,\omega_s)$ is constant. Using (i), we get
  \[
    \B(\proj\psi,\proj\psi)=\lim_{k\to\infty}\B(\proj{w_k(\omega_s)},\proj{w_k(\omega_s)})=\lim_{k\to\infty}\frac{\B(\omega_s,\omega_s)}{h(\omega_k)^2}=0.
  \]
\end{proof}

\begin{remark}
	Despite of Theorems~\ref{thm:ConvRoot} and~\ref{thm:ConvWeight}, we would like to point out that for an arbitrary space-like direction $\proj\x\in\mathbb{P}V$, the orbit $W\cdot\proj\x$ does not accumulate on~$\proj Q$ in general.  This follows from the work of Calabi and Markus~\cite{calabi_relativistic_1962}.  In the paper~\cite{chen_limit_2014}, the authors investigate in detail the space-like limit directions of Lorentzian Coxeter systems.
\end{remark}

\begin{theorem}\label{thm:limitweight}
  The set of limit weights of a Lorentzian Coxeter system $(W,S)_B$ is equal to its set of limit roots. That is, $E(\Omega)=E(\Phi)$.
\end{theorem}

\begin{proof}
  We only prove one inclusion, namely $E(\Phi)\subseteq E(\Omega)$. The proof for the other inclusion works similarly.

  Consider an injective sequence of projective roots $(\proj\gamma_k)_{k\in\mathbb{N}}$ that converges to a limit root $\proj\beta$. By Theorem \ref{thm:ConvRoot}, $\proj\beta\in\proj Q$. By passing to a subsequence, we may assume that $\gamma_k=w_k(\alpha)$ for a fixed simple root $\alpha\in\Delta$ and an injective sequence $(w_k)_{k\in\mathbb{N}}$ of elements of $W$ with increasing length. For each $k\in\mathbb{N}$, we can choose a fundamental weight $\omega_k$, such that $(w_k(\omega_k))_{k\in\mathbb{N}}$ is an injective sequence of weights.  
  This can be seen from the increasing gallery distance $d(\proj C,w_k\cdot\proj C)=\ell(w_k)$ in the Coxeter complex $\C$, which guarantees an injective sequence of vertices of $\C$, corresponding to the sequence $(w_k(\omega_k))_{k\in\mathbb{N}}$.
  By passing again to a subsequence, we may assume that $\omega_k=\omega$ for a fixed fundamental weight $\omega\in\Delta^*$, and that the sequence $(w_k(\omega))_{k\in\mathbb{N}}$ converges to a limit $\proj\psi\in\proj Q$. By Theorem \ref{thm:ConvWeight}, one has $\proj\psi\in\proj Q$. The bilinear form $\B(\alpha,\omega)$ equals $0$ or $1$, $h(w_k(\alpha))$ tends to $+\infty$ by Theorem~\ref{thm:ConvRoot1}, and $h(w_k(\omega))$ tends to $-\infty$ by Theorem~\ref{thm:ConvWeight1}. Consequently, we have:
  \[
    \B(\proj\beta,\proj\psi)=\lim_{k\to\infty}\B(\proj{w_k(\alpha)},\proj{w_k(\omega)})=\lim_{k\to\infty}\frac{\B(\alpha,\omega)}{h(w_k(\alpha))h(w_k(\omega))}=0.
  \]
  Since both $\proj\beta$ and $\proj\psi$ lie in the projective light-cone $\proj Q$, we must have $\proj\beta=\proj\psi$ by Proposition \ref{prop:twolights}. Thus, $E(\Phi)\subseteq E(\Omega)$.
\end{proof}

\begin{figure}[!ht]
  \centering
  \begin{tabular}{p{.45\textwidth}p{.45\textwidth}}
    \begin{tikzpicture}[scale=2, sommet/.style={inner sep=2pt,circle,draw=black,fill=blue,thick,anchor=base}]
      \node[anchor=south west,inner sep=0pt] at (0,0) {\includegraphics[width=.45\textwidth]{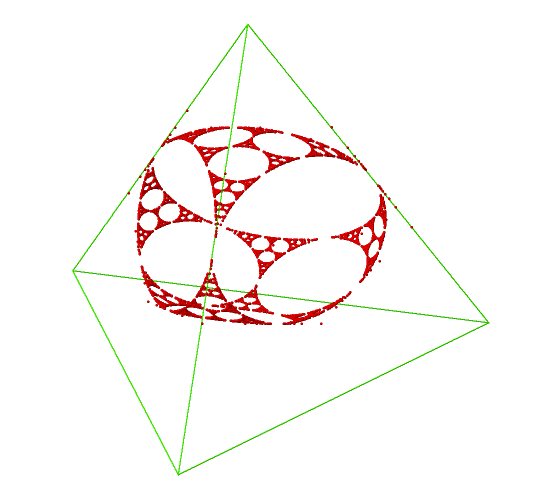}};
      \coordinate (ancre) at (0,.2\textwidth);
      \node[sommet,label=left:$s_\alpha$] (alpha) at (ancre) {};
      \node[sommet,label=right :$s_\beta$] (beta) at ($(ancre)+(0.5,0)$) {};
      \node[sommet,label=above:$s_\delta$] (delta) at ($(ancre)+(0.22,0.43)$) {} edge[thick] node[auto,swap] {$\infty$} (alpha) edge[thick] node[auto] {$\infty$} (beta);
      \node[sommet,label=below:$s_\gamma$] (gamma) at ($(ancre)+(0.125,-0.215)$) {} edge[thick] node[midway, right] {$\infty$} (delta);
    \end{tikzpicture}
    &
    \begin{tikzpicture}[scale=2]
      \node[anchor=south west,inner sep=0pt] at (0,0) {\includegraphics[width=.45\textwidth]{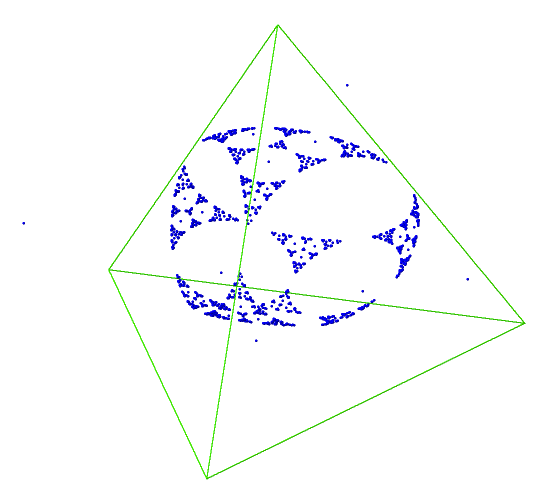}};
    \end{tikzpicture}
    \\
    \centering \small (a) Roots of depth $\leq$ 9.
    &
    \centering \small (b) Weights of ``depth'' $\leq$ 9.
  \end{tabular}
  \caption{Positive roots and space-like weights for a Lorentzian Coxeter system of rank $4$ seen in the affine space spanned by the simple roots. The Coxeter graph is shown on the upper-left corner.}
  \label{fig:limit_weights}
\end{figure}
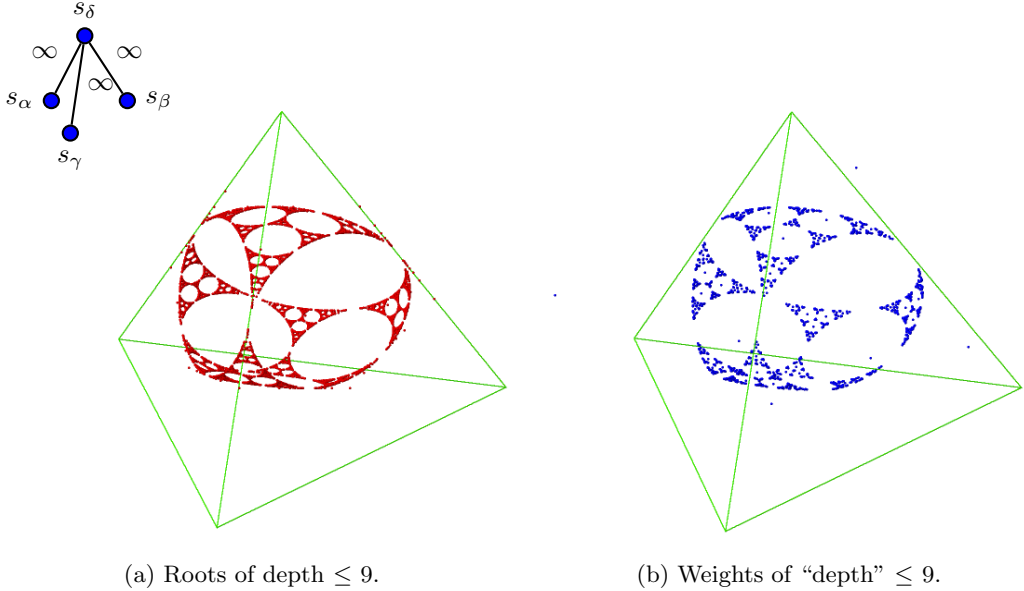

Figure~\ref{fig:limit_weights} illustrates the result of the previous theorem.
%
%
%
%
%

\begin{remark}
  It is natural to ask whether the previous result holds for non-singular geometric Coxeter systems in general. The property that distinguish Lorentz spaces is Proposition~\ref{prop:twolights}, which asserts that a totally isotropic subspace of a Lorentz space is at most of dimension~$1$. Indeed, it would be interesting to know if the equality $E(\Omega)=E(\Phi)$ holds in general. 

  For this, an answer to \cite[Question~4.9]{dyer_imaginary2_2013} would be helpful. In fact, $E(\Omega)=E(\Phi)$ holds under the assumption that $\convex(E(\Phi))\cap\proj Q=E(\Phi)$. Here is a sketch of proof: From the definition, limit weights are on the boundary of the projective Tits cone~$\proj T$; see \cite[Exercise~2.90]{abramenko_buildings_2008}. We have seen in Theorem~\ref{thm:ConvWeight2} that limit weights are on the projective isotropic cone $\proj Q$. Consequently, limit weights are in the dual of $\proj T$, which is $\convex(E(\Phi))$ by \cite[Theorem~5.1(a)]{dyer_imaginary_2013}. By assumption, we have proved that $E(\Omega)\subseteq E(\Phi)$, and the equality follows from Theorem~\ref{thm:minimal}.
\end{remark}

\subsection{Limit roots and Boyd--Maxwell ball packings}\label{ssec:Main}
We now prove Theorem~\ref{thm:main} assuming that $(W,S)_B$ is a Lorentzian Coxeter system of level~$2$.  Recall from Section~\ref{ssec:bmpacking} that~$\Omega_r$ denotes the set of space-like weights. It is the union of the orbits of space-like fundamental weights. Space-like weights in~$\Omega_r$ correspond to balls in the Boyd--Maxwell ball packing $\pack$ generated by $(W,S)_B$. From this correspondence, we prove the following theorem.
\begin{theorem}\label{thm:packing}
  The set $E(\Phi)$ of limit roots of an irreducible Lorentzian Coxeter system~$(W,S)_B$ of level~$2$ is equal to the residual set of the Boyd-Maxwell ball packing $\pack$ generated by $(W,S)_B$.
\end{theorem}

\begin{proof}
	The set $E_\Phi$ is a minimal set under the action of $W$ by Theorem~\ref{thm:minimal}.  Therefore~$E(\Phi)=E(\Omega)$ is the set of accumulation points of $\proj\Omega_r$. Since limit roots are light-like, $E(\Phi)=E(\Omega)$ is disjoint from~$\proj\Omega_r$. By Theorem~\ref{thm:maxwell}, $E(\Phi)$ is disjoint from the interiors of the balls in the Boyd--Maxwell packing~$\pack$. This proves that $E(\Phi)$ is contained in the residual set of~$\pack$.

  The other inclusion follows from the fact that $\pack$ is maximal, i.e.\ it is impossible to add any ball into the complement of $\pack$ to form a bigger packing. In other words, for a point $p$ in the residual set of $\pack$, every neighborhood of $p$ contains some ball in $\pack$. So $p$ is an accumulation point of $\proj\Omega_r$, therefore a limit root. The maximality of $\pack$ is guaranteed by~\cite[Theorem~3.3]{maxwell_sphere_1982} and \cite[Theorem~6.1]{maxwell_wythoffs_1989}.
\end{proof}

Let us now explain the relation between ball packings studied by Boyd and Maxwell and ball packings observed in the study of limit roots. Maxwell's condition of ``level~$2$'' can be interpreted as follows. Consider a Coxeter system~$(W,S)_B$ of level~$2$ with Coxeter graph $G$. Then $(W,S)_B$ is of level~$\le 2$, i.e.~removing any two vertices from $G$ leaves an affine or finite Coxeter graph. In the affine picture, this means that every $(n-2)$-face of the simplex~$\convex(\proj\Delta)$ is disjoint from, or tangent to the projective light cone $\proj Q$. Furthermore, $(W,S)_B$ is not of level $\le 1$, i.e.~there exists a vertex of $G$ whose removal does not yield an affine or finite Coxeter graph. In the affine picture, this means that some facet of the simplex $\convex(\proj\Delta)$ intersects the projective light cone~$\proj Q$ transversally.  In other words, there is at least one space-like weight.  In the point of view of \cite{hohlweg_asymptotical_2013} and \cite{dyer_imaginary2_2013}, $(W,S)_B$ is of level~$2$ if and only if $\proj Q\setminus\convex(\Delta)$ is not empty and consists of a union of disjoint open balls. Then we notice from Equation~\eqref{eqn:dualbasis} that 
\[
  \affine(\Delta\setminus\{\alpha_s\})=\affine(\Delta)\cap H_{\omega_s},\quad\forall s\in S,
\]
In other words, the supporting hyperplane $\affine(\Delta\setminus\{\alpha_s\})$ of the simplex $\convex(\Delta)$ is exactly the intersection of $\affine(\Delta)$ and the orthogonal hyperplane for the fundamental weight $\omega_s$. Therefore, the closed balls obtained by the space-like fundamental weights are exactly the closure of the open balls in $\proj Q\setminus\convex(\Delta)$. Consequently, if the irreducible Coxeter system is Lorentzian of level~$2$, the fractal structure described in Theorem \ref{thm:fractal} is the Boyd--Maxwell ball packing described in Theorem~\ref{thm:maxwell}. 

\subsection{Coxeter complex and Tangency graph}

The tangency graph $\G$ of a ball packing $\pack$ takes the balls in $\pack$ as vertices, and two vertices are connected by an edge if the corresponding balls are tangent to each other.  The tangency graph of disk packings ($2$-dimensional ball packings) is well understood, thanks to the Koebe--Andreev--Thurston's disk packing theorem, which asserts that every planar graph is the tangency graph of a disk packing.   See~\cite{stephenson_circle_2003} for a nice survey on circle packings.  However, little is known for higher dimensional ball packings; see~\cite[Section 1.3.1]{chen_thesis_2014} for a summary on previous works.  In~\cite{chen_apollonian_2013}, the first author compare the tangency graphs of Apollonian packings to $1$-skeletons of stacked polytopes, and give a forbidden subgraph characterisation for $3$-dimensional Apollonian packings.  In this part, we interpret the tangency graph of a Boyd--Maxwell ball packing in terms of the Coxeter complex of the associated Coxeter system. 

Recall that the vertices of a Coxeter complex can be colored by fundamental weights. Vertices with time- or light-like colors do not correspond to any ball in the packing, so we call them \Dfn{imaginary vertices}.  Vertices with space-like colors correspond to balls in the Boyd--Maxwell packing, so we call them \Dfn{real vertices}. For a Lorentzian Coxeter system of level 2, $\B(\omega,\omega)\le 1$ for all fundamental weights~$\omega\in\Delta^*$~\cite[Proposition 1.6]{maxwell_sphere_1982}. A vertex colored by $\omega$, such that $\B(\omega,\omega)=1$, is said to be \Dfn{surreal}, and a panel of type $\omega$ is called \Dfn{surreal}.  We use the term ``surreal'' because these vertices are not only real balls in the packing, but also guarantee tangency.  Two surreal vertices are said to be \Dfn{adjacent} if they are of the same color and belong to two adjacent chambers of the Coxeter complex $\C$ sharing a surreal panel. Then from Equation~\eqref{eqn:Bref}, we see that a pair of adjacent surreal vertices correspond to a pair of tangent balls. Finally, an edge $uv$ of type $\{\omega_u,\omega_v\}$ is called a \Dfn{real edge} if and only if $\B(\normal\omega_u,\normal\omega_v)=-1$. A real edge correspond to a pair of tangent balls in the packing.

We can now describe the tangency graph in term of the Coxeter complex.
\begin{theorem}
  \label{thm:tangency}
  Let $\pack$ be a Boyd--Maxwell ball packing generated by a Coxeter system $(W,S)_B$ of level~$2$. Let $\C$ be the Coxeter complex of $(W,S)_B$, and $\G$ be the tangency graph of $\pack$. Then the vertices of $\G$ are the real vertices of $\C$, and $uv$ is an edge of $\G$ if and only if one of the following condition is fulfilled,
  \begin{itemize}
    \item The edge $uv$ of $\C$ is real, in which case $u$ and $v$ are of different colors,
    \item The vertices $u$ and $v$ are surreal and adjacent, in which case $u$ and $v$ are of the same color.
  \end{itemize}
\end{theorem}

Edges connecting pairs of adjacent surreal vertices are present in the tangency graph $\G$ but not in the Coxeter complex $\C$, so we call them \Dfn{surreal edges}. Therefore, the tangency graph $\G$ can be constructed by taking the real vertices and real edges from the $1$-skeleton of $\C$, and add surreal edges.

\begin{proof}
  We have seen that real vertices represent balls in the Boyd--Maxwell packing, while real edges and surreal edges represent pairs of tangent balls. This was first observed by Maxwell in \cite{maxwell_sphere_1982}. It remains to prove that every pair of balls must be represented by a real edge or a surreal edge. 
  
  For two real vertices $u$ and $v$ such that $d(u,v)\ge 2$, we will prove that the balls represented by $u$ and $v$ are not tangent.  The proof is by induction on the gallery distance. The inductive step is exactly the same as in the proof of Equation (1.5) in \cite{maxwell_sphere_1982}, but we need to establish different base cases for proving strict inequalities.

  Let $u$ be a vertex of $\C$, $\omega_u\in\Delta^*$ be its color, and $s_u\in S$ and $\alpha_u\in\Delta$ be the corresponding generator and simple root, that is $s_u(\omega_u)=\omega_u-2\alpha_u$. Without loss of generality, we assume that $u=\proj\omega_u$. Let $v$ be another vertex of $\C$ and $(C_0,\dots,C_{d(u,v)})$ be a minimal gallery connecting $u$ and $v$. We may assume that $C_0$ is the fundamental chamber. Let $\omega_i$ ($1\le i\le d(u,v)$) be the type of the panel shared by $C_i$ and $C_{i-1}$, we define $w=s_1w'=s_1\dots s_{d(u,v)}$ where $s_i$ is the generator corresponding to $\omega_i$. Note that $s_1=s_u$ and $s_{d(u,v)}=s_v$. Then $C_{d(u,v)}=w\cdot C_0$, and $v=w\cdot\proj\omega_v$. Maxwell proved that~\cite[Equation~(1.6) et seq.]{maxwell_sphere_1982}
  \begin{align}
    \B(\omega_u,w(\omega_v))&=\B(\omega_u,w'(\omega_v))-2\B(\alpha_u,w'(\omega_v)), \label{eqn:induction}\\
    &\le-\sqrt{\B(\omega_u,\omega_u)\B(\omega_v,\omega_v)} \label{eqn:inequality}
  \end{align}
for any $u\ne v$. We now prove, by induction on $d(u,v)$, that the inequality \eqref{eqn:inequality} is strict for $d(u,v)\ge 2$. First, we establish the base cases.

  If $u$ and $v$ are of different colors, $\omega_u\ne\omega_v$, the base case for the induction is $d(u,v)=2$. We may assume that $w=s_us_v$, where $s_u$ and $s_v$ do not commute (otherwise $d(u,v)=0$), then
  \begin{align*}
    \B(\omega_u,w(\omega_v))&=\B(s_u(\omega_u),s_v(\omega_v))=\B(\omega_u-2\alpha_u,\omega_v-2\alpha_v),\\
    &=\B(\omega_u,\omega_v)+4\B(\alpha_u,\alpha_v)-2\B(\omega_u,\alpha_v)-2\B(\alpha_u,\omega_v).
  \end{align*}
  In the last line, the last two terms are $0$, the second term is \emph{strictly} negative since $s_u$ and $s_v$ do not commute, and the first term $\le-\sqrt{\B(\omega_u,\omega_u)\B(\omega_v,\omega_v)}$ by~\eqref{eqn:inequality}. We conclude that \[\B(\omega_u,w(\omega_v))<-\sqrt{\B(\omega_u,\omega_u)\B(\omega_v,\omega_v)}.\] Therefore, $\B(\normal\omega_u,\normal{w(\omega_v)})<-1$, and the balls represented by $u$ and $v$ are not tangent.

  If $u$ and $v$ are of the same color $\omega=\omega_u=\omega_v$, the base case for the induction is $d(u,v)=3$. Let $s=s_u=s_v$ and $\alpha=\alpha_u=\alpha_v$. We may assume that $w=ss's$, where $s\ne s'\in S$ and the order of $ss'$ is bigger than $3$ (otherwise $d(u,v)\le 1$). Then
  \begin{align*}
    \B(\omega,w(\omega))&=\B(s(\omega),s's(\omega))=\B(s(\omega),s(\omega))-2\B(s(\omega),\alpha')^2\\&=\B(\omega,\omega)-2\B(\omega-2\alpha,\alpha')^2=\B(\omega,\omega)-8\B(\alpha,\alpha')^2,
  \end{align*}
  where $\alpha'$ is the simple root corresponding to $s'$.  In the last line, the first term is $\le 1$ by \cite[Proposition 1.6]{maxwell_sphere_1982}. As for the second term, since the order of $ss'$ is bigger than~$3$, we have $\B(\alpha,\alpha')<-1/2$, so $8\B(\alpha,\alpha')^2>2$.  We conclude that $\B(\omega,w(\omega))<-1$, so $\B(\normal\omega, \normal{w(\omega)})<-1$, therefore the balls represented by $u$ and $v$ are not tangent.  
	
	In Equation \eqref{eqn:induction}, the second term $\B(\alpha_u,w'(\omega_v))$ is non-negative \cite[Corollary~1.8]{maxwell_sphere_1982}. We then use Equation \eqref{eqn:induction} for the induction, and conclude that if $d(u,v)\ge 2$, the corresponding balls are not tangent, so $uv$ does not correspond to any edge of~$\G$. Therefore, the only possible edges are the real and the surreal edges.
\end{proof}

\begin{corollary}
  For an irreducible Lorentzian Coxeter system $(W,S)_B$ of level~$2$, the projective Tits cone $\proj T\subset\mathbb{P}V$ is an edge-tangent infinite polytope. That is, every edge of $\proj T$ is tangent to the projective light-cone $\proj Q$. Furthermore, the $1$-skeleton of~$\proj T$ is the tangency graph of the Boyd--Maxwell packing~$\pack$ generated by $(W,S)_B$.
\end{corollary}

\begin{proof}
  Vertices of $\proj T$ are projective weights. No edge of $\proj T$ is disjoint from $\proj Q$, otherwise two balls in the packing $\pack$ will overlap. No edge of $\proj T$ intersect $\proj Q$ transversally because $\proj Q\subset\proj T$ by~\cite[Corollary 1.3]{maxwell_sphere_1982}. Finally, an edge of $\proj T$ that is tangent to $\proj Q$ correspond to a pair of tangent balls in~$\pack$.
\end{proof}

\subsection{Limit roots and Boyd--Maxwell ball clusters}\label{ssec:cluster}

We now finish the proof of Theorem~\ref{thm:main}.

For a Lorentzian Coxeter system of level $\le 1$, every facet of $\convex(\proj\Delta)$ is disjoint from, or tangent to $Q$. Since there are no space-like weights, the Boyd--Maxwell ball cluster is empty. Therefore $E(\Phi)=\proj Q$, as observed in \cite{hohlweg_limit_2013,dyer_imaginary2_2013}. In this case, the boundary of the Tits cone is the light cone.

For a Lorentzian Coxeter system of level~$\ge 3$, the space-like weights still represent $(n-2)$-dimensional balls, but some balls will intersect each other; see Figure~\ref{fig:patterns}(b) for an example.  To generalize Theorem~\ref{thm:packing} to Boyd--Maxwell ball clusters, most of the arguments and discussions in Section~\ref{ssec:Main} apply. However, slight modifications are necessary. 

First of all, in a Boyd--Maxwell ball cluster, we claim that no two balls intersect deeply.
\footnote{This claim is wrong. The correct claim is that there is no containment in the ball cluster.  An erratum is appended to the end of the manuscript.}
Recall that two balls intersect deeply if one is contained in the other, or if their boundary intersect at an obtuse angle, in which case the bilinear form of the corresponding space-like weights is positive. Our claim is a consequence of the following lemma, taken from the proof of Theorem~1.9 in \cite{maxwell_sphere_1982}.

\begin{lemma}[{\cite[Equation (1.5)]{maxwell_sphere_1982}}]\label{lem:falseproof}
  Let $\omega,\omega'\in \Omega$ be two distinct weights of a Lorentzian Coxeter group. Then $\B(\omega,\omega')\le 0$.
\end{lemma}

Correspondingly, we say that a ball cluster is \Dfn{maximal} if it is impossible to add any additional ball into the cluster without deeply intersecting any other ball. The maximality of Boyd--Maxwell packings is again guaranteed by the following generalized version of \cite[Theorem~3.3]{maxwell_sphere_1982}.

\begin{lemma}\label{lem:maxcluster}
  Let $(W,S)_B$ be an irreducible Lorentzian Coxeter system of level~$2$ or higher. If $\cone(\Omega)=\cone(\Omega_r)$, then the Boyd--Maxwell ball cluster generated by $(W,S)_B$ is maximal.
\end{lemma}

Maxwell's proof of \cite[Theorem~3.3]{maxwell_sphere_1982} applies directly to this generalized version, and the assumption of this lemma is verified by \cite[Theorem~6.1]{maxwell_wythoffs_1989}. All other arguments in the proof of Theorem~\ref{thm:packing} generalize directly, which completes the connection between Theorem \ref{thm:fractal} and \ref{thm:maxwell}, and the proof of Theorem~\ref{thm:main}.

%% file: List.tex
The list of level-$0$ Coxeter graphs can be found in \cite[Chapter~2]{humphreys_reflection_1992}. As observed by Maxwell \cite{maxwell_sphere_1982}, a graph of level~$2$ is either connected, or obtained by adding an isolated vertex to a graph of level $1$. Coxeter graphs of level $1$ are necessarily connected. A complete list is given by Chein in \cite{chein_recherche_1969} using a FORTRAN program; see also \cite[Section~6.9]{humphreys_reflection_1992}.

Connected Coxeter graphs of level~$2$ are manually enumerated by Maxwell in \cite{maxwell_sphere_1982}. He finds 323 Coxeter graphs of level~$2$, some of which correspond to the same packing. He then gives a list of 165 graphs, representing different packings generated by these graphs. We follow the suggestion of Maxwell and realize a computer verification of the list along the lines of \cite{chein_recherche_1969}. The list of level~2 Coxeter graphs is given in Appendix. The current section is dedicated to the description of the algorithm. The algorithm consists of two parts:
\begin{description}
  \item[\bf Nomination] A reasonably short list of candidates covering all the possible Coxeter graphs of level 2 is generated. This is to avoid checking all the graphs with less than 11 vertices.
  \item[\bf Recognition] Every nominated candidate is passed to a recognition algorithm, and is eliminated if it is not a Coxeter graph of level 2.
\end{description}

\subsection{Recognition algorithm}
The recognition algorithm is used to eliminate false candidates, and to generate the list of level-$1$ Coxeter graphs, which helps nominating candidates. Instead of following the combinatorial algorithm described in \cite{chein_recherche_1969}, our algorithm takes advantage of developments in computer science. To tell if a matrix $M$ is positive-semidefinite, we use the computer algebra system Sage to calculate the eigenvalues of $M$, and look at the sign of the smallest eigenvalue $\lambda$. If $\lambda\ge 0$, then~$M$ is positive-semidefinite. Since the considered matrices are quite small (size at most $10\times 10$), this process is done fast.

Now consider a Coxeter graph $G$ with its associated matrix $B$.  Checking if $G$ is of finite or affine type is equivalent to checking the positive-semidefiniteness of $B$ as described above. Checking if $G$ is of level $\le r$ asks to check the positive-semidefiniteness for all the $(n-r)\times(n-r)$ principle minors of $B$. Consequently, checking if $G$ is of level~$2$ requires to check if it is of level~$\le~2$ but not of level $\le 1$.

\begin{remark}
Sage can numerically calculate the eigenvalues in double precision, which gives 15-17 significant decimal digits. Since the calculation is not in arbitrary precision, it may happen that, for an eigenvalues that equals zero, the program finds a non-zero eigenvalue that is very close to $0$. This is however not a problem. In fact, for every finite or affine Coxeter graph with at most~10 vertices, the non-zero eigenvalues of its bilinear form are all bigger than $0.003$ according to our test. Therefore, double precision suffices, if all output with absolute value $<0.001$ are regarded as zero.
\end{remark}

\subsection{Nomination of candidates}

The nomination of candidates is more technical. In the spirit of~\cite{chein_recherche_1969}, we first do some graph theoretical analysis.

As observed by Maxwell~\cite{maxwell_sphere_1982}, a graph of three vertices is of level~$2$ if it contains a dotted edge, a graph of four vertices is of level~$\le 2$ if and only if it contains no dotted edge. It remains to consider graphs of five or more vertices. A level-$2$ graph with at least five vertices does not contain any dotted edge, and the only admissible labels for an edge is $3$, $4$, $5$ or $6$. 

Consider a Coxeter graphs $G$ of level~$2$. Since $G$ is not of level $\le 1$, the deletion of some vertex~$u$ from $G$ leaves a graph $G-u$ that is neither finite nor affine. Then $G-u$ is necessarily of level~$1$. As mentioned before, $G-u$ must be connected. Therefore, a Coxeter graph of level~$2$ can be obtained by connecting a vertex to a Coxeter graph of level $1$. Now let $v$ be any vertex of $G$. Every connected component of $G-v$ must be of level $\le 1$, therefore in Liste I, Liste II or the list in Appendice of \cite{chein_recherche_1969}. All these Coxeter graph of level $\le 1$ have at most one cycle, except for three graphs of level $1$, namely the complete graph~$K_4$, complete graph minus an edge $K_4-e$ and the complete bipartite graph $K_{2,3}$. 

\subsubsection{Graphs constructed from special graphs}
\begin{figure}[p]
  \centering
  \includegraphics[width=\textwidth]{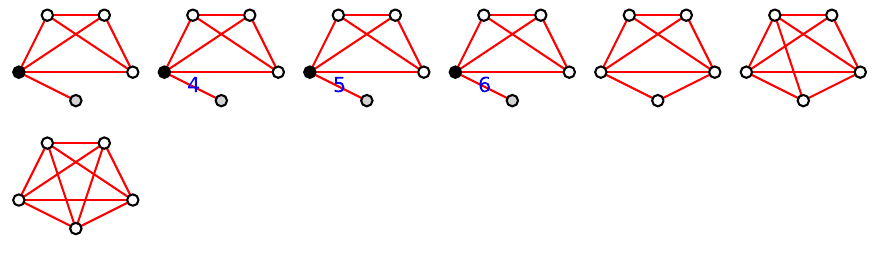}
  \caption{Graphs constructed from $K_4$.}
  \vspace{10pt}
  \includegraphics[width=\textwidth]{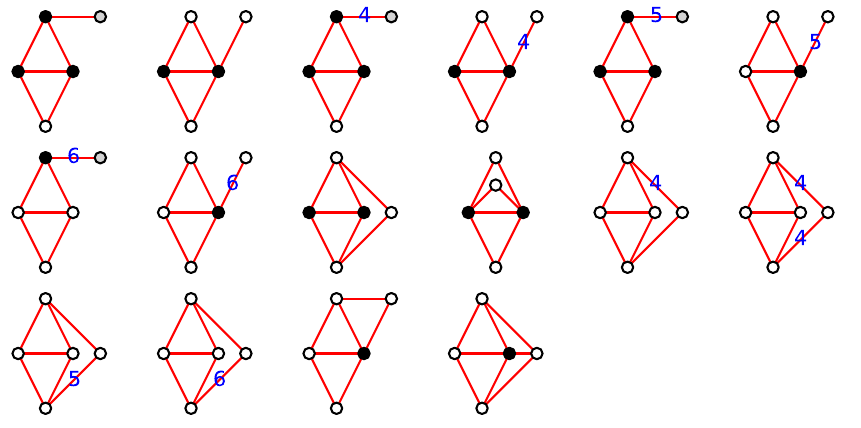}
  \caption{Graphs constructed from $K_4-e$.}
  \vspace{10pt}
  \includegraphics[width=\textwidth]{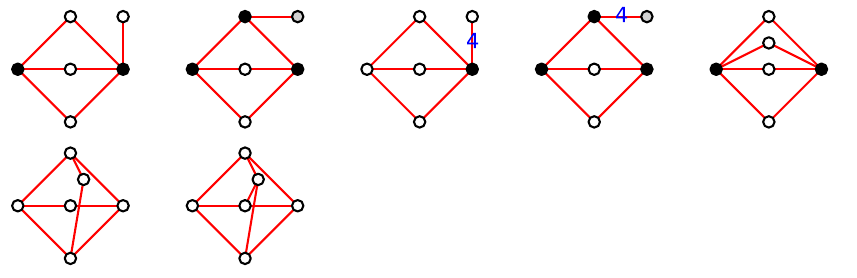}
  \caption{Graphs constructed from $K_{23}$.}
  \label{fig:S}
\end{figure}

If one component of $G-v$ contains more than one cycle, then $G$ is obtained by adding a vertex to $K_4$, $K_4-e$ or $K_{2,3}$, putting any admissible label (3, 4, 5 or 6) to the new edges. This forms our first class of candidates. After passing through the recognition algorithm, Coxeter graphs of level~$2$ constructed in this way are listed in Figure~\ref{fig:S}.
\subsubsection{Graphs with two cycles}
\begin{figure}[!ht]
  \centering
  \begin{tikzpicture}[scale=0.6]
    \draw (-4,1.2) circle (.8);
    \draw (-4,-1.2) circle (.8);
    \draw (-4,.4)--(-4,-.4);
    \draw (0,1) circle (1);
    \draw (0,-1) circle (1);
    \draw (4,0) circle (2);
    \draw (2,0)--(6,0);
  \end{tikzpicture}
  \caption{\label{fig:CCproof} The three possible forms for a level-$2$ Coxeter graph with $2$ cycles.}
\end{figure}
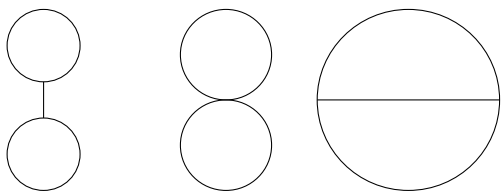

If $G-v$ contains none of the three special graphs, the argument in \cite[Section~3.2]{chein_recherche_1969} applies and we conclude the following. If $G$ has at least $5$ vertices, it has at most $2$ cycles. If the number of cycles is exactly $2$, the degree of a vertex of $G$ is at least $2$. Therefore, for a Coxeter graph with $2$ cycles, we have the three possibilities shown in Figure~\ref{fig:CCproof}.

\begin{figure}[p]
  \centering
  \includegraphics[width=.85\textwidth]{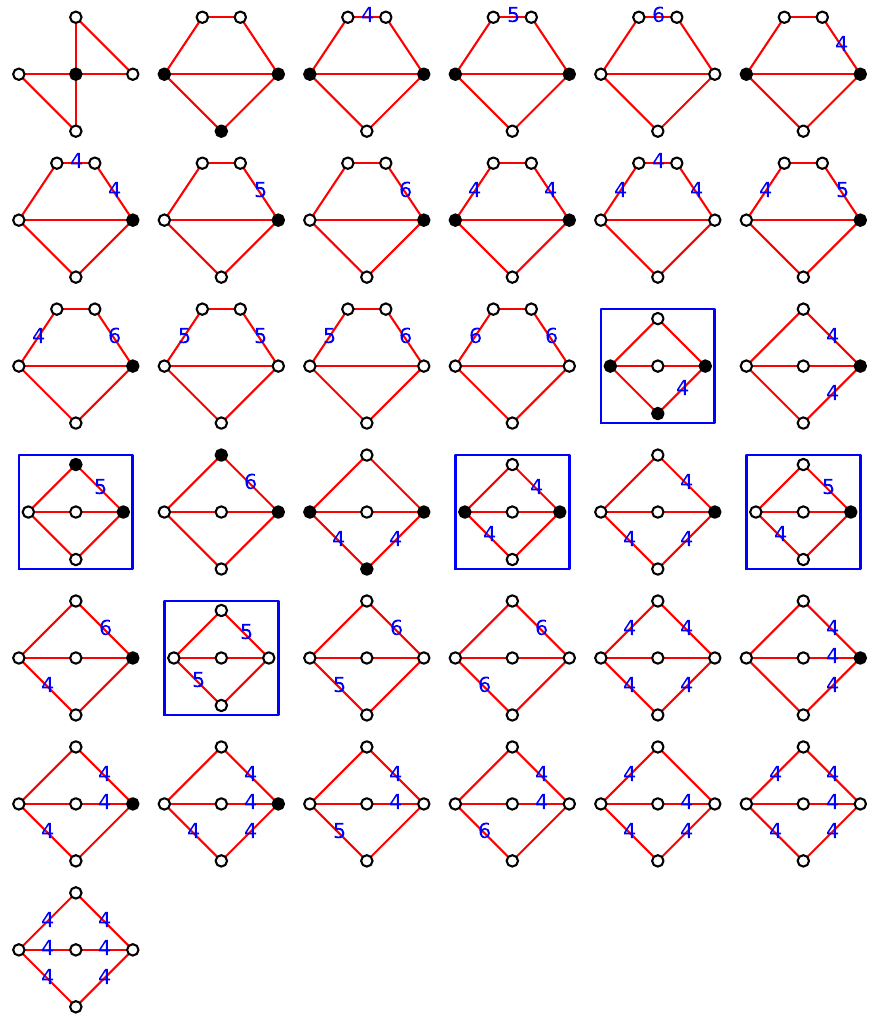}
  \includegraphics[width=.85\textwidth]{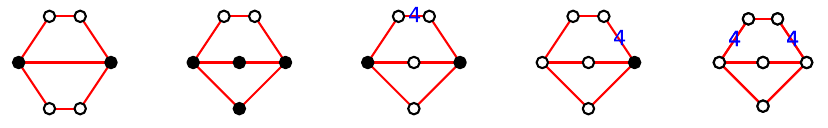}
  \includegraphics[width=.85\textwidth]{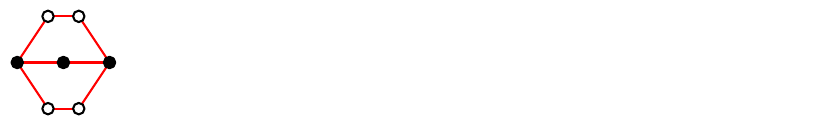}
  \includegraphics[width=.85\textwidth]{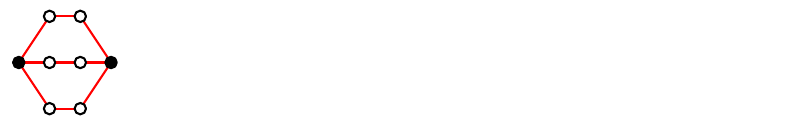}
  \caption{Graphs with two cycles.}
  \label{fig:CC}
\end{figure}

We rule out the case on the left, since deletion of two vertices from one of the two cycles leaves a graph that is not of level $0$. For the case in the middle, if any of the two cycles contains more than three vertices, deletion of two vertices on that cycle leaves a graph that is not of level $0$. The only nominated candidate is therefore the butterfly graph, i.e. two cycles of length~$3$ sharing a vertex. The butterfly graph is then confirmed by the recognition algorithm as a level-$2$ graph. For the case on the right, if any of the three paths contains more than two vertices, not counting the ends, deletion of two vertices on that path leaves a graph that is not of level $0$. Furthermore, at least two of the three paths contains at least one vertex, otherwise the graph is not simple. Graphs satisfying these two conditions, with any admissible label (3, 4, 5 or 6), are nominated as candidates. After passing through the recognition algorithm, Coxeter graphs of level~$2$ with two cycles are listed in Figure~\ref{fig:CC}.

\subsubsection{Graphs with one cycle}
A graph with only one cycle is either a cycle itself, or formed by attaching some paths to the cycle, i.e., connecting one end of the path to a vertex on the cycle. In the second case, we call the pending paths ``tails'', and the length of the tail is one plus the length of the path.

If a Coxeter graph of level $1$ has at most one cycle, there are three possibilities: a tree, a cycle, or a cycle with one tail of length $1$. If a Coxeter graph of level~$2$ has exactly one cycle, there are four possibilities: a cycle, a cycle with one tail of length $1$, a cycle with two tails of length $1$, or a cycle with one tail of length $2$. One verifies that a graph can not be of level~$2$ if it has more or longer tails.

A Coxeter graph of level~$2$ can be formed in the following ways:
\begin{enumerate}
  \item Take a tailed cycle of level $1$, then append an edge to the tail, with any admissible label (3, 4, 5 or 6). The result is a cycle with a tail of length two.  
  \item Take a tailed cycle of level $1$, then attach an edge to any vertex on the cycle, with any admissible label. The result is a cycle with two pending edges. 
  \item Take a cycle of level $1$, and attach an edge to any vertex on the cycle, with any admissible label. The result is a cycle with a tail of length $1$. 
  \item Take a tree of level $1$, then add a new vertex and connect it to any two leaves (vertices of degree $1$) of the tree, putting any admissible label to the new edges. It suffices to consider trees with three leaves, since cycles with more than two tails are not of level~$2$, and cycles with two tails are all considered in the previous case. Among cycles with one tail, we only nominate those with a tail of length one, since cycles with a tail of length two are all considered in the second case, and a longer tail length is not allowed. 
  \item Take a path of level $1$, then add a new vertex and connect it to the two ends of the path, putting any admissible label to the new edges. The result is a cycle. 
  \item Take a path of level $1$, then add a new vertex and connect it to the second and the last vertex on the path, putting any admissible label to the new edges. The result is a cycle with a tail of length~$1$. However, it turns out that all level 2 Coxeter graphs of this form have been previously nominated, and we find no new graph by this method.
\end{enumerate}

\begin{figure}[p]
  \centering
  \includegraphics[width=\textwidth]{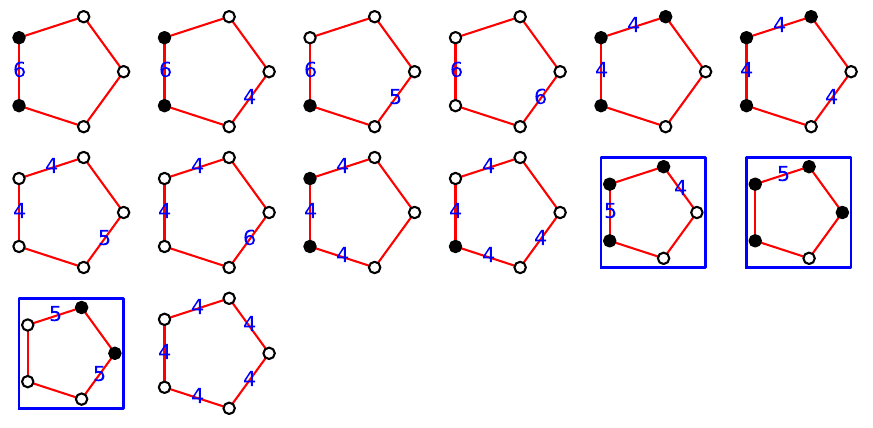}
  \includegraphics[width=\textwidth]{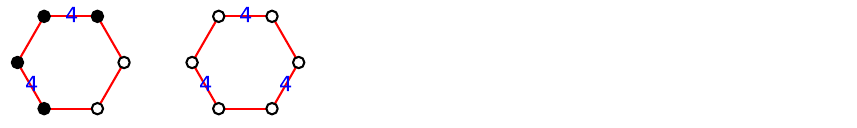}
  \includegraphics[width=\textwidth]{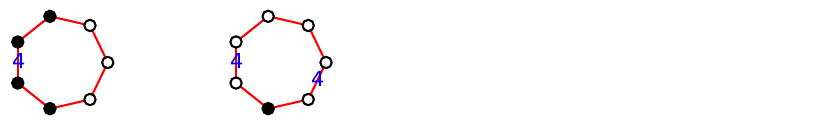}
  \caption{Cycles}
  \label{fig:C}
\end{figure}

\begin{figure}[p]
  \centering
  \includegraphics[width=\textwidth]{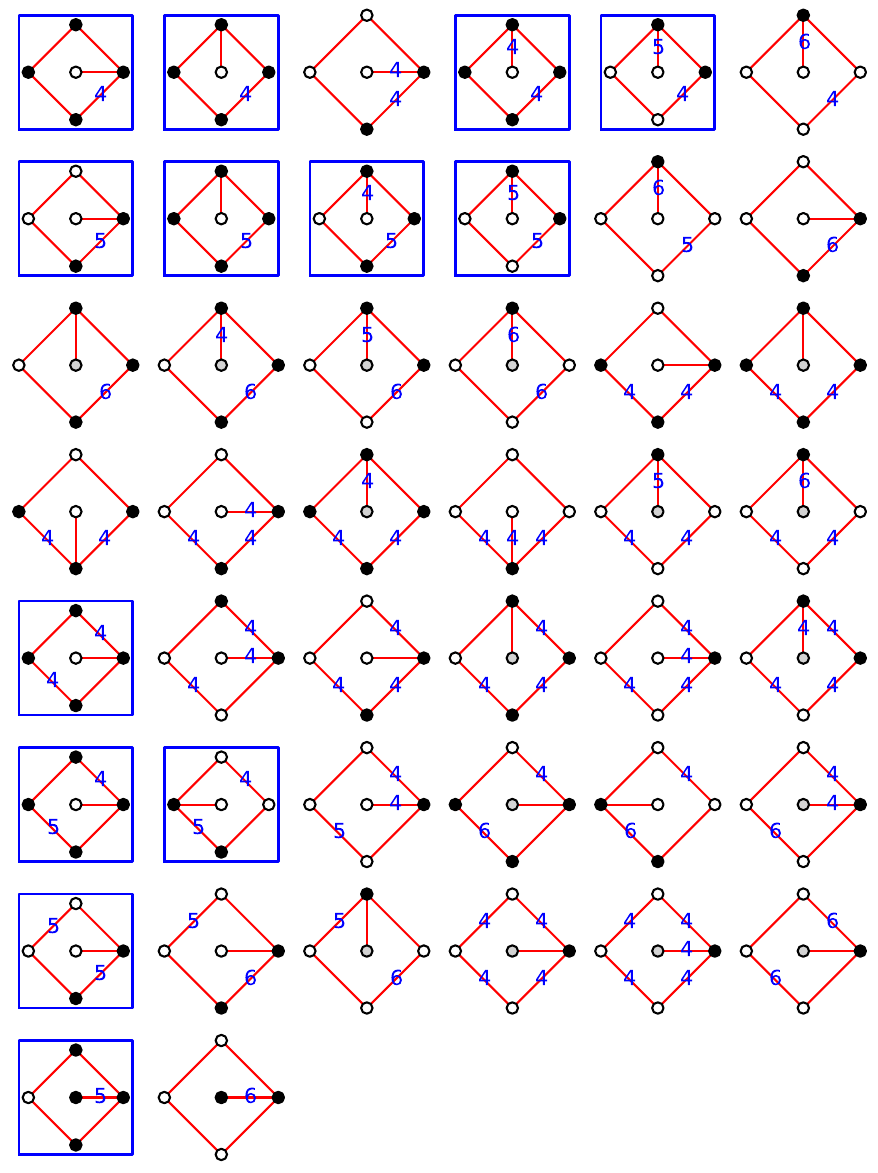}
  \caption{Cycles with one tail of length $1$ ($5$ vertices)}
  \label{fig:Ct-5}
\end{figure}

\begin{figure}[p]
  \centering
  \includegraphics[width=\textwidth]{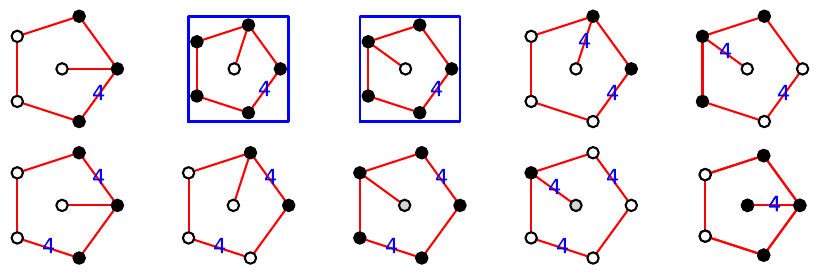}
  \includegraphics[width=\textwidth]{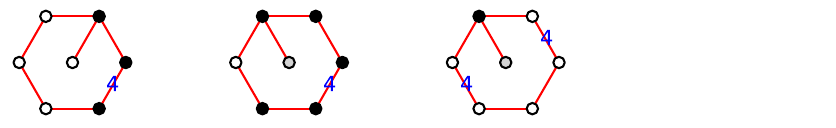}
  \includegraphics[width=\textwidth]{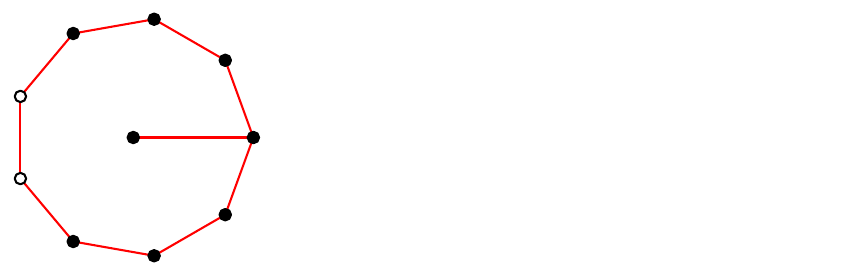}
  \caption{Cycles with one tail of length $1$ ($>5$ vertices)}
  \label{fig:Ct}
\end{figure}

\begin{figure}[p]
  \centering
  \includegraphics[width=\textwidth]{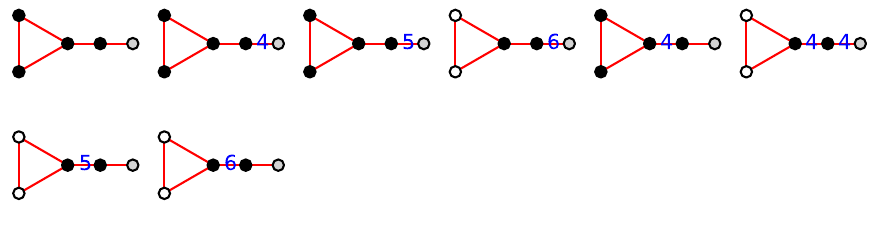}
  \includegraphics[width=\textwidth]{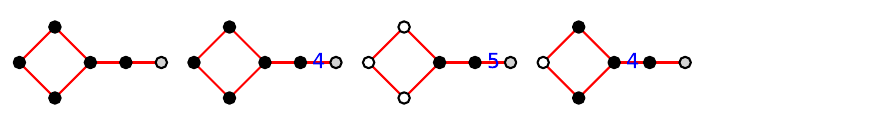}
  \includegraphics[width=\textwidth]{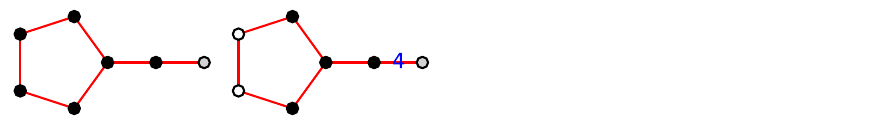}
  \includegraphics[width=\textwidth]{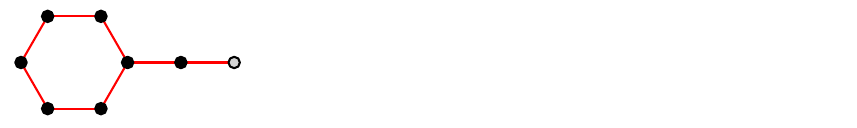}
  \includegraphics[width=\textwidth]{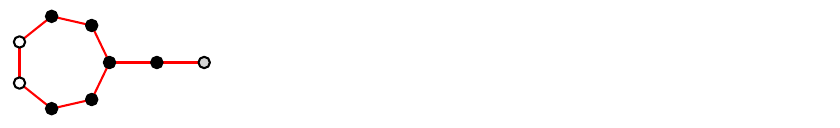}
  \caption{Cycles with one tail of length $2$.}
  \label{fig:Ct2}
\end{figure}

\begin{figure}[p]
  \centering
  \includegraphics[width=.95\textwidth]{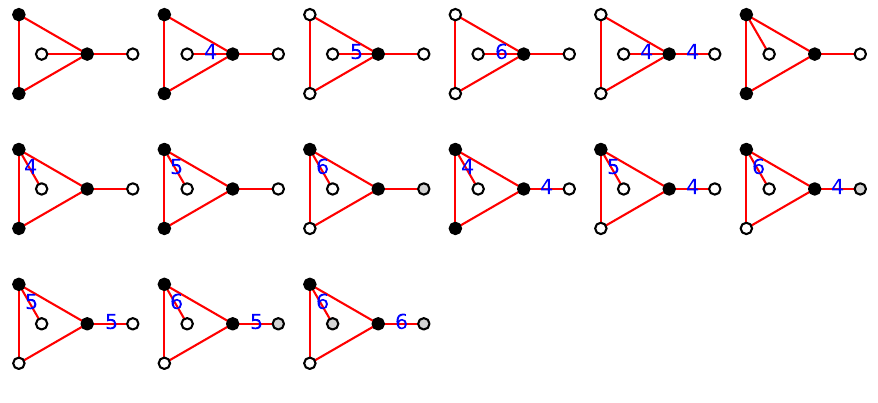}
  \includegraphics[width=.95\textwidth]{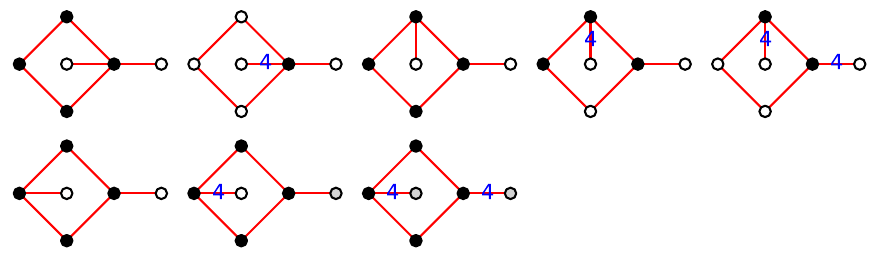}
  \includegraphics[width=.95\textwidth]{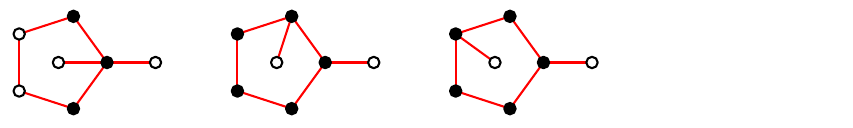}
  \includegraphics[width=.9\textwidth]{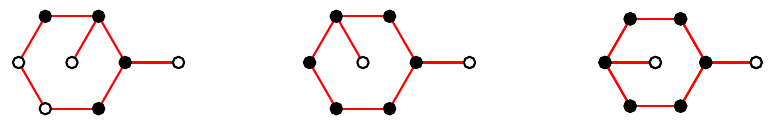}
  \includegraphics[width=.9\textwidth]{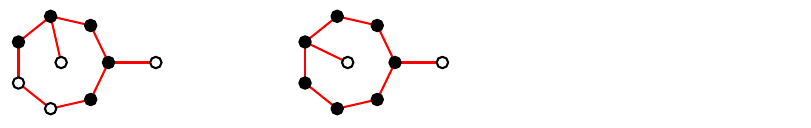}
  \includegraphics[width=.9\textwidth]{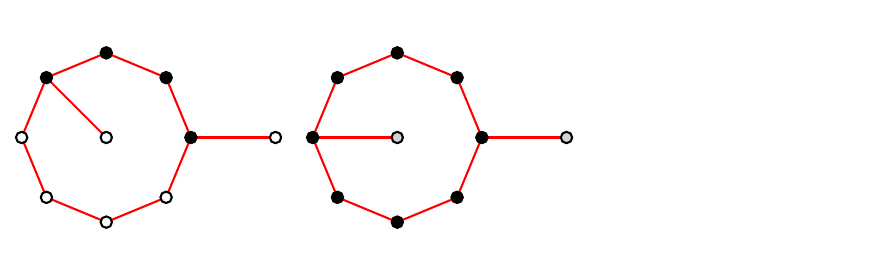}
  \caption{Cycles with two tails of length $1$.}
  \label{fig:Ctt}
\end{figure}

After passing through the recognition algorithm, Coxeter graphs of level~$2$ in form of a cycle are listed in Figure~\ref{fig:C}; those in form of a cycle with one tail of length $1$ are listed in Figure~\ref{fig:Ct-5} and~\ref{fig:Ct}; those in form of a cycle with one tail of length $2$ are listed in Figure~\ref{fig:Ct2}; those in form of a cycle with two tails of length $1$ are listed in Figure~\ref{fig:Ctt}.

\subsubsection{Graphs in form of a tree}
\begin{figure}[p]
  \centering
  \includegraphics[width=\textwidth]{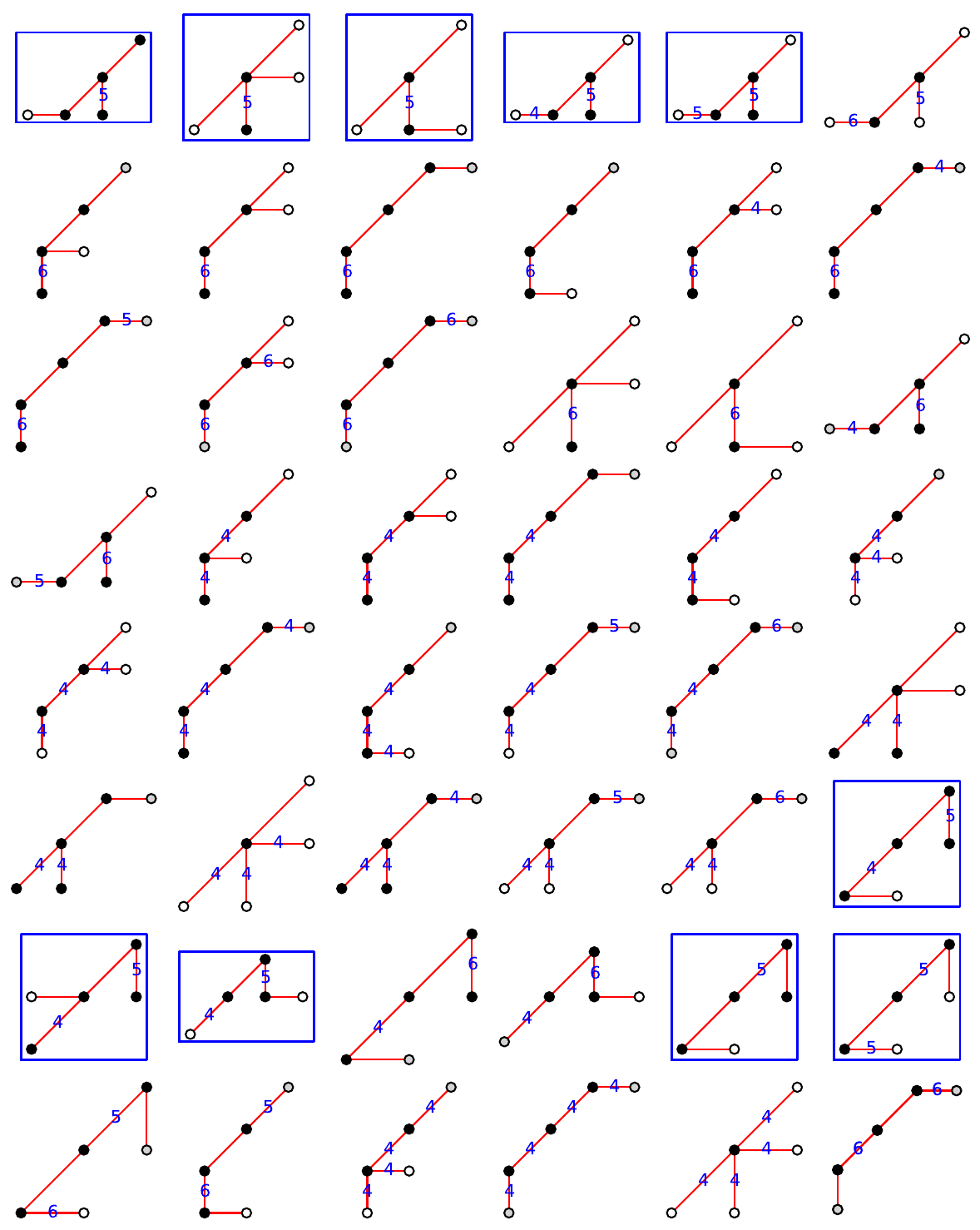}
  \caption{Trees ($5$ vertices)}
  \label{fig:T-5}
\end{figure}

\begin{figure}[p]
  \centering
  \includegraphics[width=\textwidth]{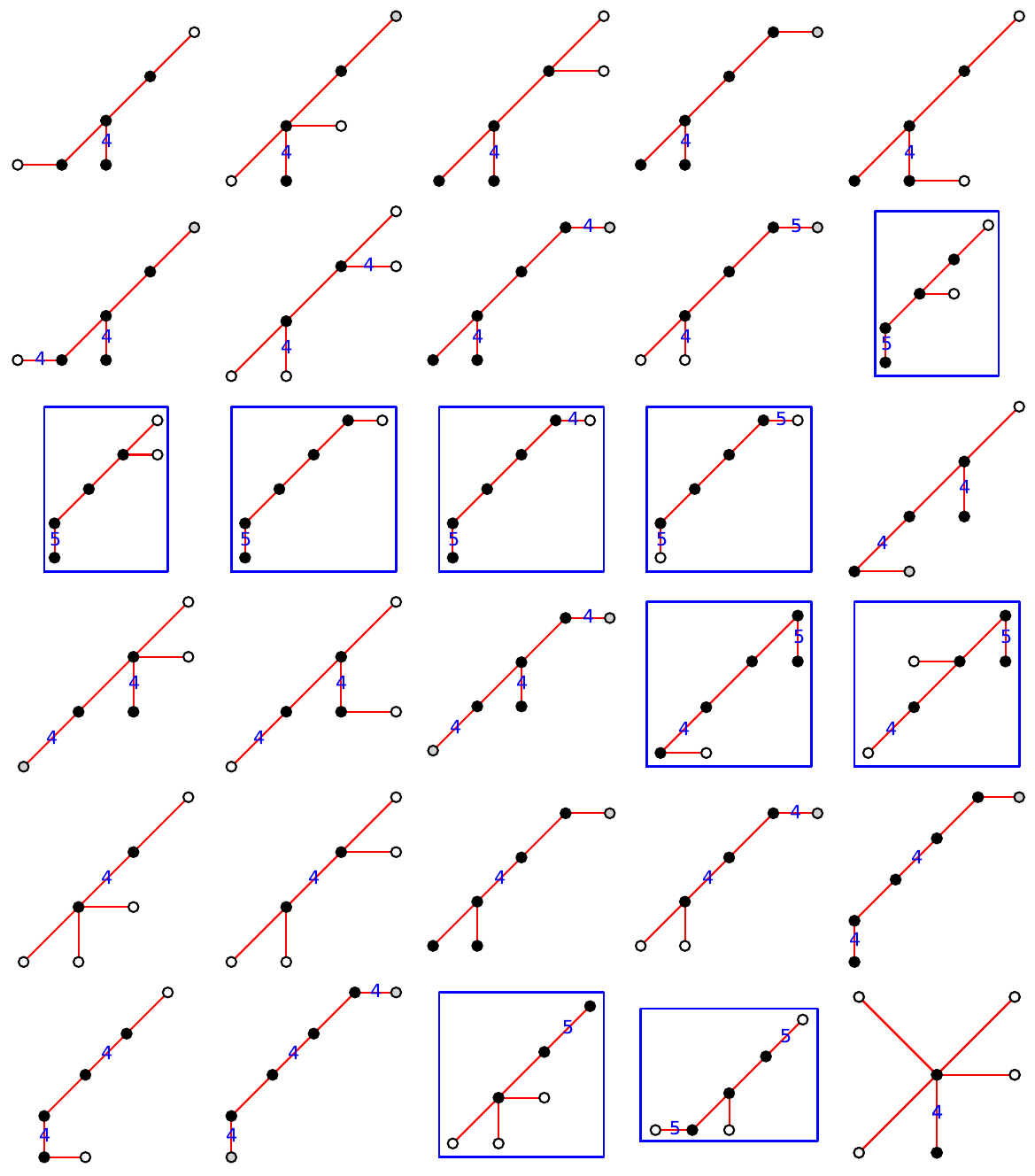}
  \caption{Trees ($6$ vertices)}
  \label{fig:T-6}
\end{figure}

\begin{figure}[p]
  \centering
  \includegraphics[width=.9\textwidth]{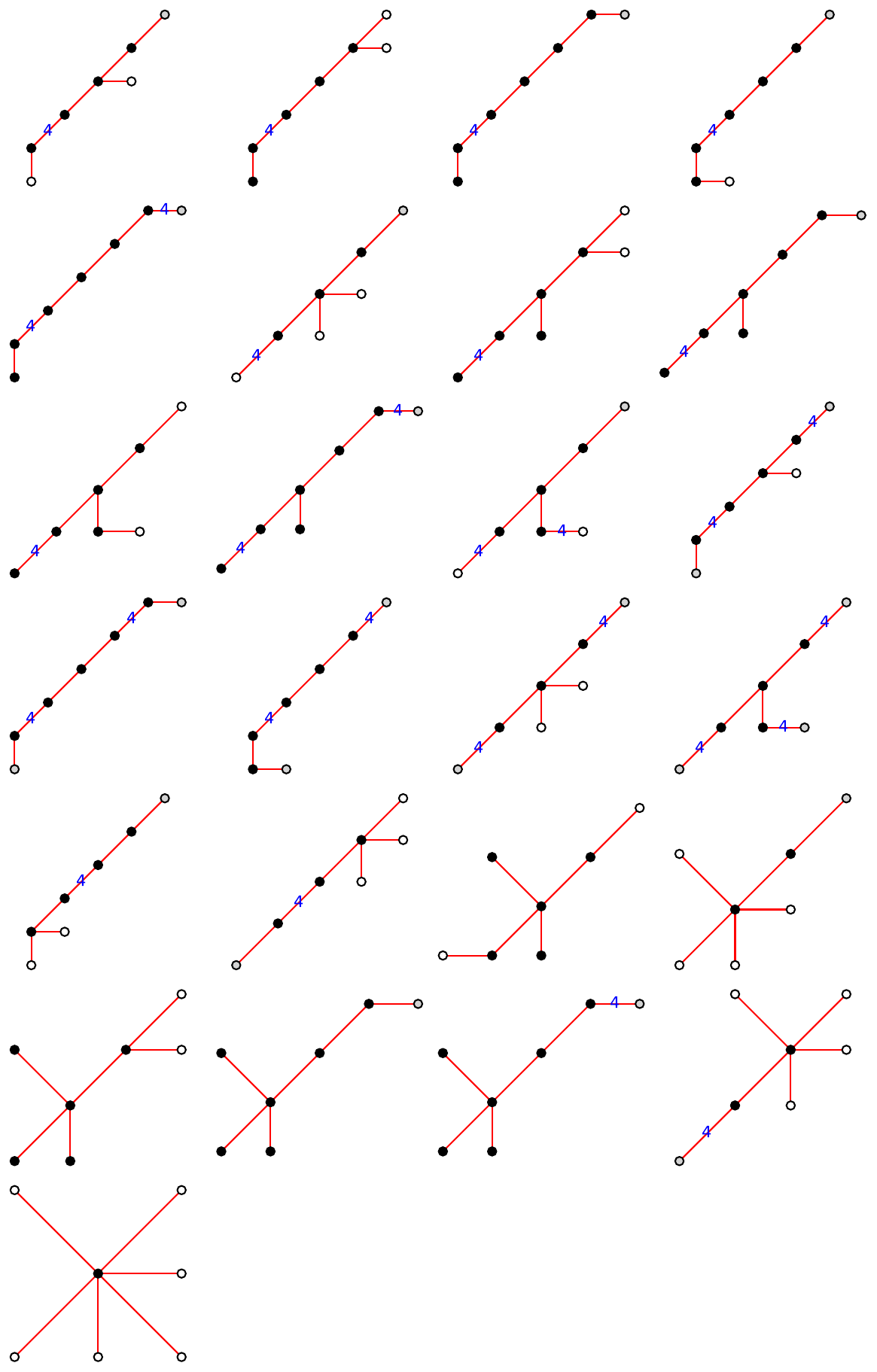}
  \caption{Trees ($7$ vertices)}
  \label{fig:T-7}
\end{figure}

\begin{figure}[p]
  \centering
  \includegraphics[width=.9\textwidth]{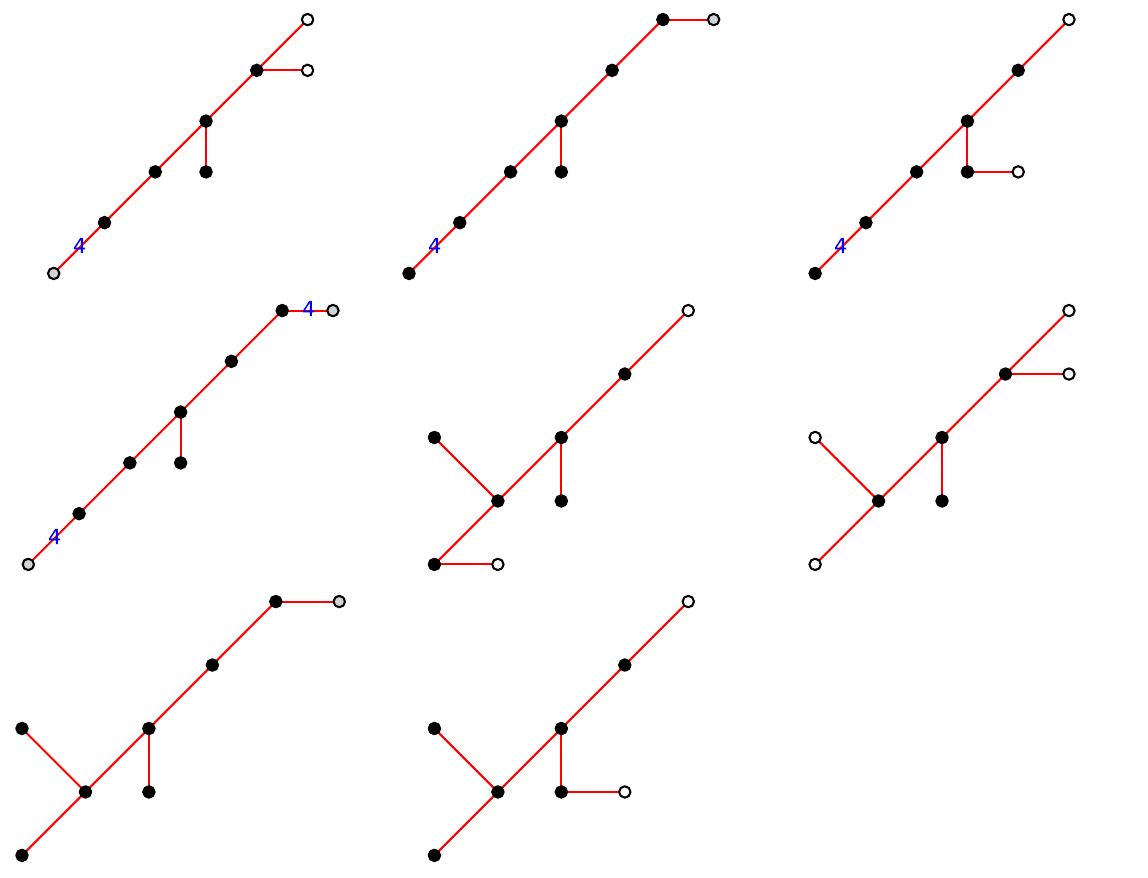}
  \includegraphics[width=.9\textwidth]{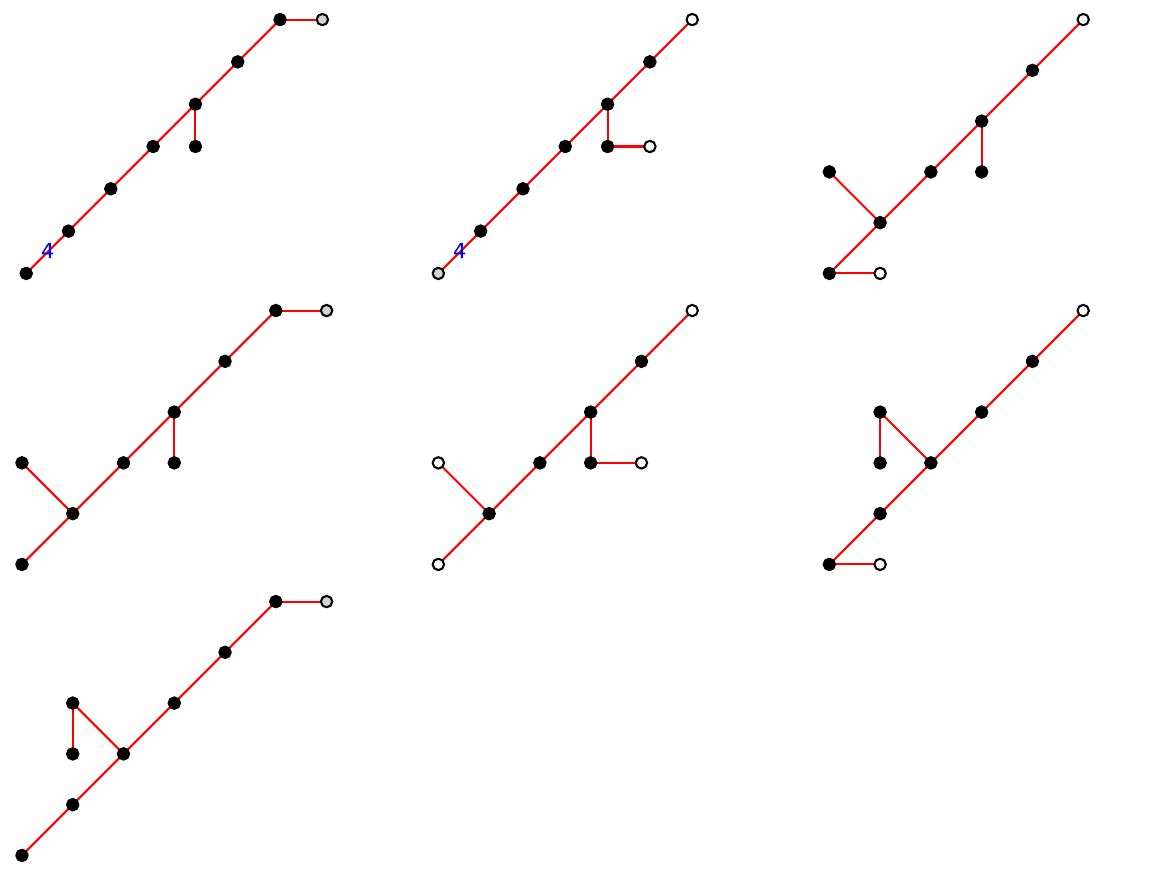}
  \caption{Trees ($8$ or $9$ vertices)}
  \label{fig:T-89}
\end{figure}

\begin{figure}[p]
  \centering
  \includegraphics[width=\textwidth]{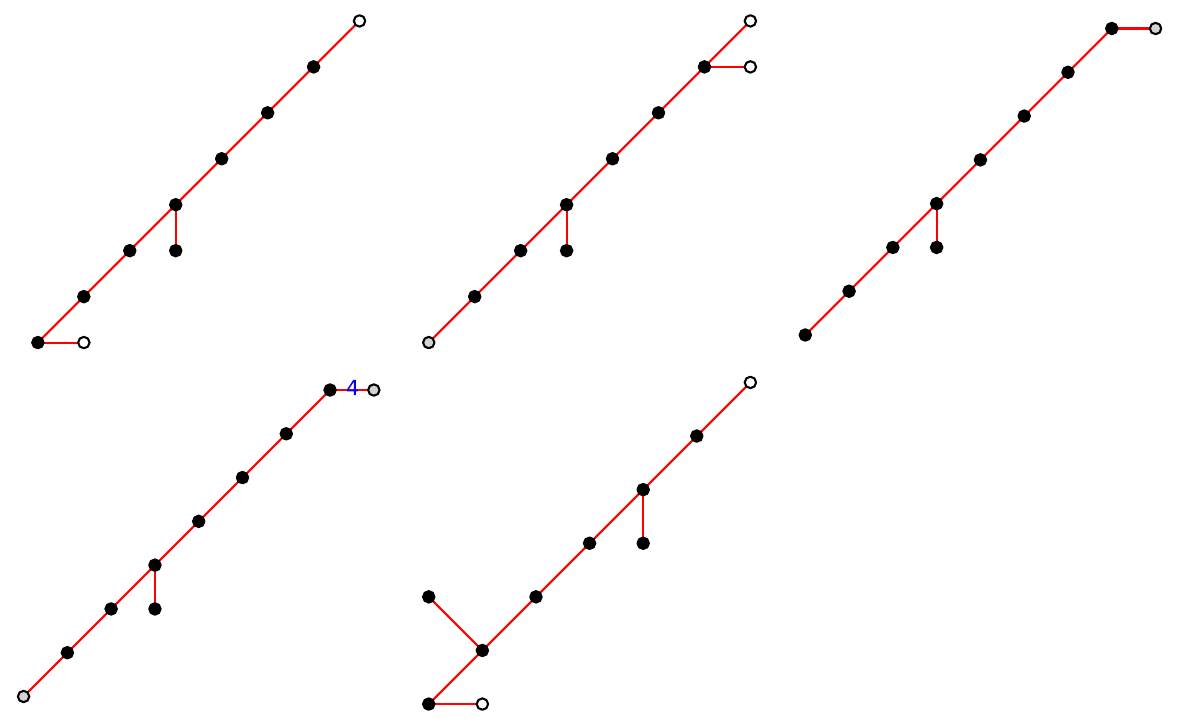}
  \includegraphics[width=\textwidth]{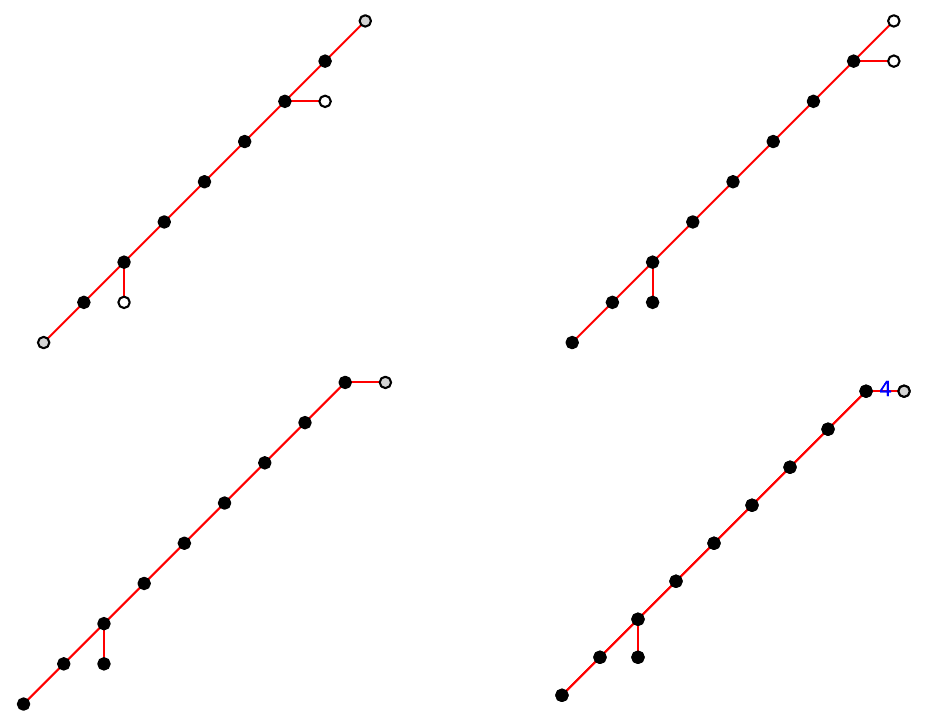}
  \caption{Trees ($10$ or $11$ vertices)}
  \label{fig:T-1011}
\end{figure}

A tree can be formed by attaching an edge to any vertex on a tree of level $1$, with any admissible label (3, 4, 5 or 6).  
After passing through the recognition algorithm, Coxeter graphs of level~$2$ in form of a tree are listed in Figures~\ref{fig:T-5} to~\ref{fig:T-1011}.

\subsection{Some remarks on the list}\label{ssec:remarks}

All Coxeter graphs of level~$2$ found by our algorithm, up to graph isomorphism, are listed in the attached tables. They are grouped according to the nomination method described above, and then subgrouped by number of vertices. The figures are generated using the graph plotting function of Sage \cite{sage}.  The implementation of the algorithm in Sage is available at~\cite{chen_sage_2014}.

The list given by Maxwell in \cite{maxwell_sphere_1982} only includes graphs corresponding to the (group theoretical) maximal elements in each family of Coxeter systems that yield the same packing. Two methods for embedding a Coxeter group as a subgroup of finite index in another Coxeter group can be found in~\cite{maxwell_euler_1998}. For checking the result, these embeddings are implemented in the program, and successfully reproduced every graph in Maxwell's list. There are 326 graphs in the present list, while Maxwell's list contains 323 Coxeter graphs of level~$2$. There are three more rank $5$ Coxeter graphs of level~$2$ in the new list. However, since Maxwell did not list all the graphs that he found, we can not specify which graphs are new.

Finite and affine Coxeter graphs with at most nine vertices are manually input into the program, then level-$1$ Coxeter graphs are enumerated following Chein's algorithm \cite{chein_recherche_1969}, from which we generate our candidates for Coxeter graphs of level~$2$. For some graphs in the list, this process can be seen from the arrangement of the edges. For example, for trees, the diagonal edges are from the original level-$0$ Coxeter graphs, vertical edges are added for constructing level-$1$ graphs, and horizontal edges are added for constructing level-$2$ graphs. For cycles with two tails, the cycles are from the original level-$0$ graphs, edges outside the cycle are added for constructing level-$1$ graphs, and edges in the cycle are added for constructing level-$2$ graphs. Colors of vertices indicate its role in the tangency graph: black vertices correspond to imaginary vertices, white and light-gray vertices correspond to real vertices, and gray vertices are surreal vertices. Some graphs are framed. These graphs of level~$2$ are \emph{strict} \cite[Section~1]{maxwell_sphere_1982}, meaning that deletion of any two vertices leaves a finite Coxeter graph. In the ball packing generated by a strict Coxeter graph of level~$2$, no two balls are tangent, i.e. the tangency graph is an empty graph. One can verify that the ball packing generated by a non-strict Coxeter graph of level~$2$ always contains a pair of tangent balls. 



%% file: Erratum.tex
In the published version of the article, there is a minor mistake in
Section~\ref{ssec:cluster}.  We fix the mistake in this appendix.

We claimed that, in a Boyd--Maxwell cluster, no two balls intersect deeply,
i.e.\ no ball contains another, and no two balls intersect at an obtuse angle.
In fact, it is possible that two balls intersect at an obtuse angle.  For an
example, consider the Coxeter graph
\tikz[scale=.8, baseline=.8em]{
  \fill (0,0) circle (.1);
  \fill (0,1) circle (.1);
  \fill (1,1) circle (.1);
  \fill (1,0) circle (.1);
	\draw (0,0) -- (1,0) -- node[right,midway,inner sep=1pt] {\small $\infty$} (1,1) -- (0,1);
	\draw[thick, dotted] (0,0) -- node[left,midway,inner sep=1pt] {\small $-1.5$} (0,1);
};
this graph is of level~$3$, and two of its space-like fundamental weights give
rise to a pair of balls intersecting at an obtuse angle.  The proof proposed in
the form of Lemma~\ref{lem:falseproof} refers
to~\cite[Equation~(1.5)]{maxwell_sphere_1982}, but the base case of the
induction in the original proof of Maxwell may fail for Coxeter groups of
level~$\ge 3$, as shown by the example above.

The correct claim is the following:
\begin{lemma*}
	In a Boyd--Maxwell ball cluster no ball is contained in another.
\end{lemma*}
\begin{proof}
	If a ball is contained in another, the corresponding space-like weight
	$\omega$ would be in the interior of the cone $\cone(\Omega_r)$.
	By~\cite[Theorem~6.1]{maxwell_wythoffs_1989}, $\cone(\Omega_r)$ is the Tits
	cone $\cone(\Omega)$.  As an interior point of the Tits cone, the stabilizer
	of $\omega$ in $W$ must be finite, so $\omega$ can not be space-like.  The
	contradiction proves the lemma.
\end{proof}

The false claim plays a minor role in the proof of the main result.  In fact,
Maxwell's proof of \cite[Theorem~3.3]{maxwell_sphere_1982} still applies to the
more general Lemma~\ref{lem:maxcluster}, whose assumption is again verified
by~\cite[Theorem~6.1]{maxwell_wythoffs_1989}.  All other arguments in the proof
of Theorem~\ref{thm:packing} remain valid, so Theorem~\ref{thm:main} is correct
as it is stated despite the mistake.